\documentclass[11pt]{amsart}

\usepackage{xypic}
\xyoption{all}

\swapnumbers
\newtheorem{thm}[equation]{Theorem}
\newtheorem{prop}[equation]{Proposition}
\newtheorem{lem}[equation]{Lemma}
\newtheorem{cor}[equation]{Corollary}

\newtheorem{conjecture}[equation]{Conjecture}
 
\theoremstyle{definition}
\newtheorem{defn}[equation]{Definition}
\newtheorem{remark}[equation]{Remark}
\newtheorem{example}[equation]{Example}

\newcommand{\nsubsection}{
\par\refstepcounter{equation}\medskip\noindent{\bf\theequation. }}

\newcommand{\diag}{\operatorname{diag}}
\newcommand{\id}{\operatorname{id}}

\newcommand{\Orb}{\operatorname{\bf Orb}}

\newcommand{\M}{\operatorname{M}}

\newcommand{\bbA}{{\mathbb A}}

\newcommand{\bbZ}{{\mathbb Z}}

\newcommand{\bbG}{{\mathbb G}}
\newcommand{\bbP}{{\mathbb P}}
\newcommand{\G}{\mathbb{G}}

\numberwithin{equation}{section}

\newcommand{\Ker}{\operatorname{Ker}}
\newcommand{\Aut}{\operatorname{Aut}}

\newcommand{\GL}{\operatorname{GL}}

\newcommand{\SO}{\operatorname{SO}}
\newcommand{\Orth}{\operatorname{O}}
\newcommand{\SL}{\operatorname{SL}}
\newcommand{\Symp}{\operatorname{Sp}}
\newcommand{\PGL}{\operatorname{PGL}}
\newcommand{\PSO}{\operatorname{PSO}}
\newcommand{\PGO}{\operatorname{PGO^+}}
\newcommand{\Pf}{\operatorname{Pf}}
\newcommand{\SPf}{\operatorname{GPf}}
\newcommand{\HW}{\operatorname{HW}}
\newcommand{\Spin}{\operatorname{Spin}}
\newcommand{\PGLn}{\PGL_n}
\newcommand{\PSymp}{\operatorname{PSp}}

\newcommand{\Stab}{\operatorname{Stab}}
\newcommand{\Sets}{\mathbf{Sets}}
\newcommand{\Fields}{\mathbf{Fields}_k}

\newcommand{\ed}{\operatorname{ed}}
\newcommand{\cd}{\operatorname{cd}}

\newcommand{\GLn}{\GL_n}
\newcommand{\lra}{\longrightarrow}
\newcommand{\Mat}{\operatorname{M}}

\newcommand{\trdeg}{\operatorname{trdeg}}

\begin{document}
\title[Canonical dimension, October 26, 04]{On the notion 
of canonical dimension for algebraic groups} 
\author[G. Berhuy and Z. Reichstein]{G. Berhuy and Z. Reichstein}
\address{Department of Mathematics, University of British Columbia,
Vancouver, BC V6T 1Z2, Canada}
\email{berhuy@math.ubc.ca}
\address{Department of Mathematics, University of British Columbia,
Vancouver, BC V6T 1Z2, Canada}
\thanks{Z. Reichstein was partially supported by an NSERC research grant}
\email{reichst@math.ubc.ca}


\begin{abstract} We define and study a numerical invariant of
an algebraic group action which we call the canonical dimension.
We then apply the resulting theory to the problem of computing
the minimal number of parameters required to define a generic
hypersurface of degree $d$ in $\bbP^{n-1}$.
\end{abstract}

\maketitle

\hfill{\em To M. Artin on his 70th birthday}

\tableofcontents

\section{Introduction}
Many important objects in algebra can be parametrized
by a non-abelian cohomology set of the form $H^1(K, G)$,
where $K$ is a field
and $G$ is a linear algebraic group defined over $K$.
For example, elements of $H^1(K, \Orth_n)$ can be identified with
isomorphism classes of $n$-dimensional quadratic forms over $K$, 
elements of $H^1(K, \PGL_n)$ with isomorphism classes 
of central simple algebras of degree $n$, elements of
$H^1(K, G_2)$ with isomorphism classes of octonion algebras, etc.; 
cf.~\cite{serre-gc}
or~\cite{boi}.  Recall that $H^1(K, G)$ has a  marked (split)
element but usually no group structure.
Thus, a priori there are only two types
of elements in $H^1(K, G)$, split and non-split.
However, it is often intuitively clear that
some non-split elements are closer to being split
than others. This intuitive notion can be quantified
by considering degrees or Galois groups of splitting field
extensions $L/K$ for $\alpha$;
see, e.g.,~\cite{tits},~\cite{ry2}.
Another ``measure" of how far $\alpha$ is from being split
is its essential dimension (here and in the sequel
we assume that $K$ is a finitely generated extension of an
algebraically closed base field of characteristic zero,
and $G$ is defined over $k$); for
details and further references, see Section~\ref{prel.ed}
and the first paragraph of Section~\ref{sect.functor}.

In this paper we introduce and study yet another
numerical invariant that ``measures" how far $\alpha$
is from being split. We call this new invariant the {\em canonical
dimension} and denote it by $\cd(\alpha)$.  We give several
equivalent descriptions of $\cd(\alpha)$; one of them is
that $\cd(\alpha) = \min \, \trdeg_K (L)$, where the minimum is
taken over all generic splitting fields $L/K$ for $\alpha$ (see
Section~\ref{sect.gen-spl}).  Generic splitting fields
have been the object of much research in the context of
central simple algebras (i.e., for $G = \PGLn$;
see, e.g.,~\cite{amitsur}, \cite{artin},
\cite{roquette1}, \cite{roquette2}) and
quadratic forms (i.e., for $G = \Orth_n$ or $\SO_n$;
see, e.g.,~\cite{knebusch1},~\cite{knebusch2},~\cite{ks}); 
related results for Jordan pairs can be found 
in~\cite{petersson}.
Kersten and Rehmann~\cite{kr}, who, following on the work 
of Knebusch, studied generic splitting fields 
in a setting rather similar to ours
(cf.~Remark~\ref{rem.kr}), remarked, on p. 61, that
the question of determining the minimal possible transcendence
degree of a generic splitting field (or $\cd(\alpha)$,
in our language) appears to be difficult in general.
Much of this paper may be viewed as an attempt to address
this question from a geometric point of view.

Recall that we are assuming $k$ to be an algebraically
closed base field of characteristic zero, and $K/k$ to be a
finitely generated field extension.
In this context every $\alpha \in H^1(K, G)$
is represented by a (unique, up to birational isomorphism)
generically free $G$-variety $X$, with $k(X)^G = K$;
see e.g.,~\cite[(1.3.3)]{popov}.  We will often work
with $X$, rather than $\alpha$,
writing $\cd(X, G)$ instead of $\cd(\alpha)$ and using
the language of invariant theory, rather than Galois cohomology.
An advantage of this approach is that $\cd(X, G)$ is
well defined for $G$-varieties $X$ that are not necessarily
generically free (see Definition~\ref{def1}), and
the interplay between generically free and non-generically free
varieties can sometimes be used to gain insight into their
canonical dimensions; cf., e.g., Lemma~\ref{lem5.2}.
If $S$ is the stabilizer in general position
for a $G$-variety $X$, then $\cd(X, G)$
can be related to the essential dimension of $S$.
This connection is explored in
Sections~\ref{sect.lower}--\ref{sect.rel-sect}.

In Sections~\ref{sect.defn-grp}--\ref{sect.cd1} we study
canonical dimensions of generically free $G$-varieties
or, equivalently, of classes $\alpha \in H^1(K, G)$.
We will be particularly interested in the maximal possible value
of $\cd(\alpha)$ for a given group $G$; we call this number
{\em the canonical dimension of $G$} and denote it by $\cd(G)$.
The canonical dimension $\cd(G)$, like the essential dimension $\ed(G)$,
is a numerical invariant of $G$; if $G$ is connected, both measure,
in different ways, how far $G$ is from being ``special"
(for the definition and a brief discussion of special groups,
see Section~\ref{prel.special} below).
While $\cd(G)$ and $\ed(G)$ share some common properties
(note, in particular, the similarity between the results
of Section~\ref{sect.defn-grp} in this paper and
those of~\cite[Sections 3.1, 3.2]{reichstein}), their 
numerical values do not appear to be related to each other.
For example, since $\cd(G) = 0$ for every finite group $G$ 
(see Lemma~\ref{lem6.1}(b)), the
rich theory of essential dimension for finite groups
(see~\cite{br1},~\cite{br2},~\cite[Section 8]{jly})
has no counterpart in the setting of canonical
dimension. On the other hand, our classification
of simple groups of canonical dimension $1$
in Section~\ref{sect.cd1} 
has no counterpart in the context of essential dimension,
because connected groups of essential dimension $1$ do not
exist; see~\cite[Corollary 5.7]{reichstein}.

In Section~\ref{sect.spl-fld} we prove a strong necessary
condition for $\alpha \in H^1(K, G)$ to be of canonical dimension 
$\le 2$.  A key ingredient in our proof is the 
Enriques-Manin-Iskovskih classification of minimal 
models for rational surfaces; see the proof of 
Proposition~\ref{prop.spl-deg}. In Sections~\ref{sect.a_{n-1}} 
and~\ref{sect.son} we study canonical dimensions of the groups 
$\GL_n/\mu_d$, $\SL_n/\mu_e$, $\SO_n$ and $\Spin_n$. Our arguments 
there heavily rely on the recent results of Karpenko and
Merkurjev~\cite{karpenko1}, \cite{km}, \cite{merkurjev2}.

Our definition of canonical dimension naturally extends to
the setting of functors $\mathcal{F}$ from
the category of field extensions of $k$ to the category of
pointed sets; $\cd(G)$ is then a special case of $\cd(\mathcal{F})$,
with $\mathcal{F}= H^1(_-, G)$
(see Section~\ref{sect.cd-functor}). A similar notion
in the context of essential dimension is due to
Merkurjev~\cite{merkurjev};  see also~\cite{bf2} and
the beginning of Section~\ref{sect.functor}.

In the Sections~\ref{sect.functor} -- \ref{sect.hom2}  
we apply our results 
on canonical dimension to the problem of computing
the minimal number $\ed \, [H_{n, d}]$ of independent 
parameters, required to define the general degree $d$
hypersurface in $\bbP^{n-1}$. (For a precise statement of
the problem, see Section~\ref{sect.functor}.)  We show
that if $d \ge 3$ and $(n, d) \ne (2, 3), (2, 4)$ or $(3,3)$,
our problem reduces to that of computing
the canonical dimension of the group
$\SL_n/\mu_{\gcd(n, \, d)}$. In particular, combining
Theorem~\ref{thm.hyper} with Corollary~\ref{cor6.3}, 
we obtain following theorem.

\begin{thm} \label{thm.intro} 
Let $n$ and $d$ be positive integers such that
$d \ge 3$ and $(n, d) \ne (2, 3)$, $(2, 4)$ or $(3, 3)$.
Suppose $\gcd(n, d)$ is a prime power $p^j$ for some $j \ge 0$.
Then
\[ \ed(H_{n, d}) =  
 \begin{pmatrix} n + d -1 \\ d \end{pmatrix}
- n^2 + \begin{cases} 
\text{$0$, if $j = 0$} \\
\text{$p^i - 1$, if $j \ge 1$.} 
\end{cases} \, , \]
where $p^i$ is the highest power of $p$ dividing $n$.
\end{thm}

If $d \le 2$ or $(n, d) = (2, 3), (2, 4), (3, 3)$, then
our problem reduces to computing canonical dimensions for certain
group actions that are not generically free; this is done
in Section~\ref{sect.hom2}.  Related results for $(n, d) = (2, 3)$
and $(3, 3)$ can be found in~\cite{bf1}.

\section*{Acknowledgements} We are grateful to J.-L. Colliot-Th\'el\`ene,
S. Garibaldi, D. W. Hoffmann, N. Karpenko, B. Kunyavskii, A. Merkurjev,
A. Qu\'eguiner-Mathieu, 
and J-P. Serre for helpful comments.

\section{Notation and preliminaries}

Throughout this paper we will work over an algebraically closed base
field $k$ of characteristic zero. Unless otherwise specified,
all algebraic varieties, algebraic groups, group
actions, fields and all maps between them are assumed to be
defined over $k$, all algebraic groups are assumed
to be linear (but not necessarily connected),
and all fields are assumed to be finitely generated over $k$.

By a $G$-variety we shall mean an algebraic variety $X$ with a (regular)
action of an algebraic group $G$. We will usually assume that $X$
is irreducible and focus on properties of $X$ that are preserved
by ($G$-equivariant) birational isomorphisms. In particular, we will
call a subgroup $S \subset G$ {\em a stabilizer in general position}
for $X$ if $\Stab(x)$ is conjugate to $S$ for $x \in X$
in general position; cf.~\cite[Section 7]{pv}.
As usual, if $S = \{ 1 \}$, i.e., $G$ acts freely 
on a dense open subset of $X$,
then we will say that the $G$-variety $X$ (or equivalently,
the $G$-action on $X$) is {\em generically free}.

\nsubsection {\bf Essential dimension.}
\label{prel.ed}
Let $X$ be a generically free $G$-variety. The essential dimension
$\ed(X, G)$ of $X$ is the minimal value of $\dim(Y) - \dim(G)$, where
the minimum is taken over all dominant rational maps $X \dasharrow Y$ of
$G$-varieties with $Y$ generically free. For a given algebraic group $G$,
$\ed(X, G)$ attains its maximal value in the case where $X = V$ is
a (generically free) linear representation of $G$. This value is
called the essential dimension of $G$ and is denoted by $\ed(G)$
(it is independent of the choice of $V$). For details, see
\cite[Section 3]{reichstein}.

\nsubsection {\bf Rational quotients.}
\label{prel.rat-quot}
A rational quotient for a $G$-variety $X$ is an algebraic variety $Y$
such that $k(Y) = k(X)^G$. The inclusion $k(Y) \hookrightarrow k(X)$
then induces a rational quotient map $\pi \colon X \dasharrow Y$.
Note that $Y$ and $\pi$ are only defined up to birational isomorphism;
one usually writes $X/G$ in place of $Y$.  We shall say that
{\em $G$-orbits in $X$ are separated by regular invariants}
if $\pi$ is a regular map and $\pi^{-1}(y)$ is a single $G$-orbit
for every $k$-point $y \in Y$.  By a theorem
of Rosenlicht, $X$ has a $G$-invariant dense open subset
$U$, where $G$-orbits are separated by regular invariants.
For a detailed discussion of the rational quotient
and Rosenlicht's theorem, see \cite[Section 2.4]{pv}.

\nsubsection {\bf Generically free actions and Galois cohomology.}
\label{prel.gen-free}
Let $X$ be a generically free variety. Then $X$ may be viewed as 
a torsor over the generic point of $X/G$ via the rational quotient map
$X \dasharrow X/G$. Let $\alpha$ be the class of this torsor 
in $H^1(K, G)$, where $K = k(X)^G$. This class is
explicitly constructed in~\cite[(1.3.1)]{popov}. Moreover, every
$\alpha \in H^1(K, G)$ can be obtained in this way, and the $G$-variety
$X$ can be uniquely reconstructed from $\alpha$, up to a ($G$-equivariant)
birational isomorphism; see~\cite[(1.3.2) and (1.3.3)]{popov}. 
In the sequel we shall say that $\alpha \in H^1(K, G)$ 
{\em represents} the generically free $G$-variety $X$.
 
\nsubsection {\bf Split generically free varieties.}
\label{prel.split}
Let $X$ be a generically free $G$-variety, where the $G$-orbits
are separated by regular invariants.  We will call a rational
map $s \colon X/G \dasharrow X$ {\em a rational section}
for $\pi$ if $s \circ \pi = \id$ on $X/G$.
(Note that since the fibers of $\pi$ are precisely the
$G$-orbits in $X$, $G \cdot s(X/G)$ is dense in $X$.
Consequently, some translate of $s$ will ``survive"
if $X$ is replaced by a birationally equivalent $G$-variety.)
We shall say that $X$ is {\em split} if one of the following equivalent
conditions holds:

\smallskip
(i) $X$ is birationally isomorphic to $G \times X/G$,

\smallskip
(ii) $\pi$ has a rational section,

\smallskip
(iii) $X$ represents the trivial class in $H^1(K, G)$,

\smallskip
(iv) $\ed(X, G) = 0$.

\smallskip
\noindent
For a proof of equivalence of these four conditions, see~\cite[(1.4.1)]{popov}
and~\cite[Lemma 5.2]{reichstein}. 

\nsubsection {\bf The groups $\GL_n/\mu_d$ and $\SL_n/\mu_e$.}
\label{prel.slnd}
In this section we will review known results about
the Galois cohomology
sets $H^1(K, G)$, where $G = \GL_n/\mu_d$ or $\SL_n/\mu_e$,
$\mu_d$ is the unique central cyclic subgroup of $\GL_n$ of order $d$,
and $e$ divides $n$. 

\begin{lem} \label{lem.slnd}
Let $G = \GL_n/\mu_d$ (respectively, $G = \SL_n/\mu_e$),
$f \colon G \lra \PGLn$ be the canonical projection, 
and $K/k$ be a field extension. Then

\smallskip
(a) The map $f_{\ast} \colon H^1(K, G)\to H^1(K,\PGL_n)$
has trivial kernel.

\smallskip
(b) The image of $f_{\ast}$ consists of those classes
which represent central simple algebras of degree $n$
and exponent dividing $d$ (respectively, dividing $e$).
\end{lem}

Lemma~\ref{lem.slnd} can be deduced from~\cite[Theorem 3.2]{saltman1};
for the sake of completeness, we supply a direct proof below.

\begin{proof}
(a) The exact sequence
$1 \lra \Ker(f) \stackrel{i}{\hookrightarrow}
G \stackrel{f}{\lra} \PGL_n \lra 1$
of algebraic groups, gives rise to an exact sequence
\[ H^1(K, \Ker(f)) \stackrel{i_*}{\lra}
H^1(K, G) \stackrel{f_*}{\lra} H^1(K, \PGL_n) \] 
of pointed sets; cf.~\cite[pp. 123 - 126]{serre-lf}.
It is thus enough to show that $i_*$ is the trivial map
(i.e., its image is $\{ 1 \}$).
If $G = \GL_n/\mu_d$ this is an immediate consequence 
of the fact that $\Ker(f) = \bbG_m/\mu_d$ is isomorphic 
to $\bbG_m$ and thus $H^1(K, \Ker(f)) = \{ 1 \}$.
If $G = \SL_n/\mu_e$ then $\Ker(f) = \mu_{\frac{n}{e}}$, 
and the commutative diagram
\begin{equation} \label{e.commdiagr}
\xymatrix{
1\ar@{->}[r]&\mu_n\ar@{->}[r]\ar@{->}[d]^{\times \, e}&\SL_n\ar@{->}[r]\ar@{->}[d]&\PGL_n\ar@{=}[d]\ar@{->}[r]&1\cr
1\ar@{->}[r]&\mu_{\frac{n}{e}} \ar@{->}[r]^{i}&\SL_n/\mu_e\ar@{->}[r]&\PGL_n\ar@{->}[r]&1,}
\end{equation}
of group homomorphisms induces the commutative diagram
$$\xymatrix{
H^1(K, \mu_n) \ar@{->}[r]\ar@{->}[d] & H^1(K, \SL_n)
\ar@{->}[d] \cr
H^1(K, \mu_{\frac{n}{e}}) \ar@{->}[r]^{i_*}& H^1(K, \SL_n/\mu_e) }$$
of maps of pointed sets. Since the left vertical map
is surjective (it is the natural projection
$K^{\times} /(K^{\times})^n \lra K^{\times} /(K^{\times})^\frac{n}{e}$), and
$H^1(K, \SL_n) = \{ 1 \}$ (see~\cite[p. 151]{serre-lf}),
we see that the image of $i_*$ is trivial,
as claimed.

\smallskip
(b) We will assume $G = \SL_n/\mu_e$; the case
$G = \GL_n/\mu_d$ is similar and will be left to the reader.
We now focus on the connecting maps
$$\xymatrix{
 & H^1(K, \PGLn) \ar@{->}[r]^{\delta}\ar@{=}[d] & H^2(K, \mu_n)
\ar@{->}[d]^{\times \, e} \cr
H^1(K, \SL_n/\mu_e) \ar@{->}[r]^{f_*} &
 H^1(K, \PGLn) \ar@{->}[r]^{\delta'} \ar@{->}[r] & H^2(K, \mu_{\frac{n}{e}})}
$$
induced by the diagram~\eqref{e.commdiagr}. It is well known that
$H^2(K, \mu_n)$ is the $n$-torsion part
of the Brauer group of $K$ and $\delta$ sends a central 
simple algebra $A$
to its Brauer class $[A]$; see~\cite[Section X.5]{serre-lf}.
Hence, $\delta'(A) = e \cdot [A]$,
and $\operatorname{Im}(f_*) = \Ker(\delta')$
consists of algebras $A$ of degree $n$ and exponent dividing $e$, 
as claimed.
\end{proof}

\nsubsection {\bf Special groups.}
\label{prel.special}
An algebraic group $G$ is called {\em special} if $H^1(E, G) = \{ 1 \}$
for every field extension $E/k$. Equivalently, $G$ is special 
if every generically free
$G$-variety is split. Special groups were introduced by Serre~\cite{serre1}
and classified by Grothendieck~\cite[Theorem 3]{grothendieck}
as follows: $G$ is special if and only if its maximal semisimple
subgroup is a direct product of simply connected groups
of type $\SL$ or $\Symp$; cf. also~\cite[Theorem 2.8]{pv}.
The following lemma can be easily deduced from Grothendieck's 
classification; we will instead give a proof based on Lemma~\ref{lem.slnd}.

\begin{lem} \label{lem.e}
$\GL_n/\mu_d$ is special if and only if $\gcd(n, d) = 1$.
\end{lem}

\begin{proof} If $n$ and $d$ are
relatively prime then every central simple algebra
of degree $d$ and exponent dividing $n$ is split.
By Lemma~\ref{lem.slnd}, $f_{\ast}$ has trivial image
and trivial kernel,
showing that $H^1(E, \GL_n/\mu_d) = \{ 1 \}$ for every $E$,
i.e., $\GL_n/\mu_d$ is special.
Conversely, suppose $e = \gcd(n, d) > 1$.
Let $E = k(a, b)$,
where $a$ and $b$ are algebraically independent variables
over $k$, and $D = (a, b)_e$ = generic symbol algebra of degree $e$.
Then $A = M_{\frac{n}{e}}(D)$ is a central simple algebra 
of degree $n$ and exponent $e$, with center $E$. 
This algebra defines a class in $H^1(E, \PGL_n)$; since $e$ divides $d$,
Lemma~\ref{lem.slnd} tells us that this class is the image
of some $\alpha \in H^1(E, \GL_n/\mu_d)$. Since $A$ is not split,
$\alpha \ne 1$, and hence $\GL_n/\mu_d$ is not special, as claimed.
\end{proof}

\section{The canonical dimension of a $G$-variety}
\label{sect.def}

\begin{defn} \label{def.fandF}
Let $X$ be an irreducible $G$-variety (not necessarily generically free).
We shall say that a rational map $F \colon X \dasharrow X$
is a {\em canonical form} map if $F(x) = f(x)\cdot x$ for
some rational map $f \colon X \dasharrow G$. Here
we think of $F(x)$ as a ``canonical form" of $x$.
Note that $F$ and $f$ will usually not be $G$-equivariant.
\end{defn}

\begin{remark} \label{rem1}
If the $G$-action on $X$ is generically free and
$F \colon X \dasharrow X$ is a rational map then
the following conditions are equivalent:

\smallskip
(a) $F$ is a canonical form map,

\smallskip
(b) $F(x) \in G \cdot x$ for $x \in X$ in general position,

\smallskip
(c) $\pi \circ F = \pi$, where $\pi \colon X \dasharrow X/G$ is
the rational quotient map.

\smallskip
\noindent
The equivalence of (a) and (b) follows from the fact that 
the rational quotient map $\pi \colon X \dasharrow X/G$ 
is a $G$-torsor over a dense open subset of $X/G$; 
cf. Section~\ref{prel.gen-free}. The equivalence of (b) and (c)
is a consequence of the theorem of Rosenlicht mentioned in
Section~\ref{prel.rat-quot}.
\end{remark}

\begin{remark} \label{rem2}
If the $G$-action on $X$ is generically free then
the argument we used to prove that 
(a) $\Leftrightarrow$ (b) in Remark~\ref{rem1}
also shows that the rational map
$f \colon X \dasharrow G$ in Definition~\ref{def.fandF}
is uniquely determined by $F$.  On the other hand,
if the $G$-action on $X$ is not generically free
then this may no longer be the case.  For example, if the $G$-action 
on $X$ is trivial then every $f \colon X \dasharrow G$ gives
rise to the trivial canonical form map $F = \id_X \colon X \dasharrow X$.
This means that in working with canonical form maps, we will 
want to keep track of both $F$ and $f$. In the sequel we will
often say that $f \colon X \dasharrow G$ {\em induces} 
$F \colon X \dasharrow X$ if $F(x) = f(x) \cdot x$ for $x \in X$ 
in general position.
\end{remark}

\begin{example} \label{ex0.5}
Let $X = \Mat_n$, with the conjugation
action of $G = \GLn$. We claim that the rational map
$F \colon \Mat_n \dasharrow \Mat_n$, taking $A$ to its
companion matrix
\[ F(A) = \begin{pmatrix}
0 & 0 & \dots & 0 & -c_n \\
1 & 0 & \dots & 0 & -c_{n-1} \\
0 & 1 & \dots & 0 & -c_{n-2} \\
\vdots & \vdots & \vdots & \vdots & \vdots \\
0 & 0 & \dots & 1 & - c_1 \end{pmatrix} \, , \]
is a canonical form map. Here
$t^n + c_1 t^{n-1} + \dots + c_n = \det(tI - A)$
is the characteristic polynomial of $A$.

To prove the claim, fix a non-zero column vector $v \in k^n$
and define $f(A)$ as the matrix whose columns are $v, Av, \dots, A^{n-1}v$.
It is now easy to see that
$f \colon A \mapsto f(A)$ is a rational map
$\Mat_n \dasharrow \GL_n$, and $f(A) \cdot A = f(A)^{-1} A f(A)$
is the companion matrix $F(A)$.
\qed
\end{example}

Our definition of a canonical form map is quite
general; for example it includes the trivial case,
where $f(x) = 1_G$ and thus $F(x) = x$ for every $x \in X$.
Usually we would like to choose $f$ so that
the canonical form of every element lies in some
subvariety of $X$ of small dimension. With this in mind,
we give the following:

\begin{defn} \label{def1}
The canonical dimension $\cd(X, G)$ of a $G$-variety $X$ is
defined as
\[ \text{$\cd(X, G) =$ min $\{ \dim \, F(X)  - \dim(X/G) \}$,} \]
where the minimum is taken over all canonical form
maps $F \colon X \dasharrow X$. If the $G$-action on $X$ 
is generically free, and $X$ represents $\alpha \in H^1(K, G)$
(see Section~\ref{prel.gen-free}), we will 
also write $\cd(\alpha)$ in place of $\cd(X, G)$.
\end{defn}

Note that the symbol $\cd$ does not stand for and should not
be confused with cohomological dimension.

\begin{lem} \label{lem2}
The integer $\cd(X, G)$ is the minimal value of
$\dim \, F(G \cdot x)$ for $x \in X$ in general position.
Here the minimum is taken over all canonical form
maps $F \colon X \dasharrow X$.
\end{lem}

\begin{proof} Let
$\pi \colon X \dasharrow X/G$ be the rational quotient
map for the $G$-action on $X$. Then
for any canonical form map $F \colon X \dasharrow X$, we have
$\pi = \pi \circ F$.  In particular,
$\pi \, F(X)$ is dense in $X/G$. 
Applying the fiber dimension theorem to
\[ \pi_{| \, F(X)} \colon F(X) \dasharrow X/G \, , \]
we see that
\[ \dim \, F(X)  - \dim \, X/G  = 
\dim \, \overline{F(X)} \cap G \cdot x =
\dim \, F(G \cdot x) \]
for $x \in X$ in general position. By Definition~\ref{def1},
$\cd(X, G)$ is the minimal value of this quantity, as $F$ ranges
over all canonical form maps $F \colon X \dasharrow X$.
\end{proof}

\begin{example} \label{ex1} Let $X = \Mat_n$, with the conjugation
action of $G = \GLn$. Then the canonical form map $F$
constructed in Example~\ref{ex0.5} takes every orbit in $X$
to a single point. This shows that $\cd(\Mat_n, \GL_n) = 0$.
The same argument shows that $\cd(\Mat_n, \PGL_n) = 0$;
cf. also Lemma~\ref{lem3} below.
\end{example}

\section{First properties}

\nsubsection {\bf Subgroups.}

\begin{lem} \label{lem2.1}
If $X$ is a $G$-variety and $H$ is a closed subgroup of $G$
then
\[ \cd(X, G) + \dim \, X/G \le \cd(X, H) + \dim \, X/H \, . \]
\end{lem}

\begin{proof} The left hand side is the minimal value of $\dim \, F(X)$,
as $F$ ranges over canonical form maps $F \colon X \dasharrow X$
induced by $f \colon X \dasharrow G$. The right hand side is
the same, except that $f$ is only allowed to range over rational maps
$X \dasharrow H$. Since there are more rational maps from $X$ to $G$
than from $X$ to $H$, the inequality follows.
\end{proof}

\nsubsection {\bf Connected components.}

\begin{lem} \label{lem2.15} Let $X$ be a $G$-variety and let
$G^0$ be the connected component of $G$. Then
$\cd(X, G^0) = \cd(X, G)$.
\end{lem}

\begin{proof}
The inequality $\cd(X, G) \le \cd(X, G^0)$ follows from
Lemma~\ref{lem2.1}, with $H = G^0$.
To prove the opposite inequality, let
$F \colon X \dasharrow X$ be a canonical
form map such that $\dim \, F(G \cdot x) = \cd(X, G)$ for $x \in X$
in general position.  Suppose $F$ is induced by
a rational map $f \colon X \dasharrow G$, as in Definition~\ref{def.fandF}.
Since $X$ is
irreducible, the image of $f$ lies in some irreducible component
of $G$. Let $g$ be an element of this component.
Then we can replace $f$ by $f' \colon X \dasharrow G^0$, where
$f'(x) = g^{-1} f(x)$, and $F$ by $F' \colon X \dasharrow X$
given by $F'(x) = f'(x) \cdot x = g^{-1} \cdot F(x)$. (Note that here $g$
is independent of $x$.) Since $F'(G \cdot x)$ is a translate of
$F(G \cdot x)$, we conclude that
$\cd(X, G^0) \le \dim \, F'(G^0 \cdot x) \le \dim \, F'(G \cdot x) =
\dim \, F(G \cdot x) = \cd(X, G)$.
\end{proof}

\nsubsection {\bf Direct products.}

\begin{lem} \label{lem.product}
Let $X_i$ be a $G_i$-variety for $i=1, 2$, $G = G_1 \times G_2$ and
$X = X_1 \times X_2$.  Then $\cd(X ,G) \le \cd(X_1,G_1)+\cd(X_2,G_2)$.
\end{lem}

\begin{proof} If $F_i \colon X_i \dasharrow X_i$ are canonical form maps
induced by $f_i \colon X_i \dasharrow G_i$ (for $i = 1, 2$) then
$F = (F_1, F_2) \colon X \dasharrow X$ is a canonical form map
induced by
$f = (f_1, f_2) \colon X = X_1 \times X_2 \dasharrow G_1 \times G_2$.
Clearly,
\[ F(G \cdot x) = F_1(G_1 \cdot x_1) \times F_2(G_2 \cdot x_2)  \]
for any $x = (x_1, x_2)$ and thus
$\dim(F \cdot x) = \dim(F_1 \cdot x_1) + \dim(F_2 \cdot x_2)$.
The desired inequality now follows from Lemma~\ref{lem2}.
\end{proof}

\nsubsection {\bf Split varieties.}

\begin{lem} \label{lem2.2}
Let $X$ be a generically free $G$-variety and let
$\pi \colon X \dasharrow X/G$ be the rational quotient map.

\smallskip
(a) If $X$ is split (cf. Section~\ref{prel.split}) then $\cd(X, G) = 0$.

\smallskip
(b) Suppose $G$ is connected. Then the converse to part (a) holds as well.
\end{lem}

\begin{proof}
(a) Since $X$ is split, we may assume $X = G \times X_0$, where $X_0 = X/G$;
see Section~\ref{prel.split}(i). The map $F \colon X \lra X$, given by
$F \colon (g, x_0) \mapsto (1_G, x_0)$ is clearly a canonical form map
(see Remark~\ref{rem1}), with $\dim \, F(X) = \dim(X/G)$, 
and the desired equality follows.

\smallskip
(b) After replacing $X$ be a $G$-invariant dense open subset, we may
assume that the $G$-orbits in $X$ are separated by regular invariants.
Suppose $\cd(X, G) = 0$, i.e., $\dim \, F(X) = \dim(X/G)$
for some canonical form map $F \colon X \dasharrow X$.
It is enough to show that $\pi_{|F(X)} \colon F(X) \dasharrow X/G$
is a birational isomorphism. Indeed, if we can prove this then
\[ \pi_{|F(X)}^{-1} \colon
X/G \stackrel{\simeq}{\dasharrow} F(X) \hookrightarrow X \]
will be a rational section (as defined in Section~\ref{prel.split}).

To prove that $\pi_{|F(X)}$ is a birational isomorphism, consider the
fibers of this map. If $x \in X$ is a point in general position
and $y = \pi(x) \in X/G$ then  $\pi_{|F(X)}^{-1}(y) = F(G \cdot x)$.
Since $G$ is connected, $G \cdot x$ is irreducible, and so is
$F(G \cdot x)$.  On the other hand, since $\cd(X, G) = 0$,
$\pi_{|F(X)}^{-1}(y)$ is 0-dimensional. We thus conclude
that $\pi_{|F(X)}^{-1}(y)$ is a single $k$-point for $y \in X/G$
in general position. Hence, $\pi_{|F(X)}$ is a birational
isomorphism (cf., e.g.,~\cite[Section I.4.6]{humphreys}),
and the proof is complete.
\end{proof}

\nsubsection {\bf Normal subgroups}

\begin{lem} \label{lem3} Let $\alpha \colon G \lra \overline{G}$
be a surjective map of algebraic groups and $H = \Ker(\alpha)$.
Suppose $X$ is a $\overline{G}$-variety or, equivalently,
a $G$-variety with $H$ acting trivially.  Then

\smallskip
(a) $\cd(X, G) \ge \cd(X, \overline{G})$.

\smallskip
(b) If $H$ is special then $\cd(X, G) = \cd(X, \overline{G})$.
\end{lem}

\begin{proof} Part (a) follows from Definition~\ref{def1}, because
$f \colon X \dasharrow G$ and
\[ \overline{f} = f \negthickspace \negthickspace \pmod{H}
\colon X \dasharrow \overline{G} \]
give rise to the same canonical form map $F \colon X \dasharrow X$.

(b) Reversing the argument of part (a), it suffices to show that
every rational map $\overline{f} \colon X \dasharrow \overline{G}$
can be lifted to $f \colon X \dasharrow G$.

Since $\alpha \colon G \lra \overline{G}$ separates 
the orbits for the right $H$-action on $G$, it is 
a rational quotient map for this action: cf, \cite[Lemma 2.1]{pv}.
If $H$ is special then $\alpha$ has a rational section $\beta \colon
\overline{G} \dasharrow G$. (Note that $\beta$ is a rational map of
varieties but not necessarily a group homomorphism.)
Moreover, for any $g_0 \in G$, the map 
$\beta_{g_0} \colon \overline{G} \dasharrow G$,
given by
$\beta_{g_0} \colon \overline{g} \mapsto g_0^{-1} 
\beta(\alpha(g_0) \overline{g})$
is also a rational section of $\alpha$. 
After replacing $\beta$ by $\beta_{g_0}$,
for a suitable $g_0 \in G$, we may assume that
$\overline{f}(X)$ does not lie entirely in
the indeterminacy locus of $\beta$.
Now $f = \beta \circ \overline{f} \colon X \dasharrow G$ is 
the desired lifting of $\overline{f} \colon X \dasharrow \overline{G}$.
\end{proof}

\begin{prop} \label{prop.product2}
Let $X$ be a $G$-variety and $H$ be a closed normal
subgroup of $G$.  If $H$ is special and 
the (restricted) $H$-action on $X$ is generically 
free then $\cd(X, G) = \cd(X/H, G) = \cd(X/H, G/H)$.
\end{prop}

Strictly speaking, the rational quotient variety
$X/H$ is only defined up 
to birational isomorphism. However, there exists 
a birational model $Y$ of $X/H$, such that the $G$-action
on $X$ descends to a (regular) $G/H$-action (or equivalently,
a regular $G$-action) on $Y$;
see~\cite[Proposition 2.6 and Corollary to Theorem 1.1]{pv}. 
The symbol $X/H$ in the statement of 
Proposition~\ref{prop.product2} denotes $Y$ as 
above; any two such models will are birationally 
isomorphic as $G$-varieties. 

\begin{proof} 
The equality $\cd(X/H, G) = \cd (X/H, G/H)$ follows
from Lemma~\ref{lem.product}(b); we shall thus focus on
proving that $\cd(X, G) = \cd(X/H, G)$.

After replacing $X$ by an $G$-invariant
open subset, we may assume that the quotient map
$\pi \colon X \lra X/H$ is regular. If 
$\pi' \colon X/H \dasharrow Z$ is
a rational quotient map for the $G$-action on $X/H$
then $\pi' \circ \pi \colon X \dasharrow Z$ is a 
rational quotient map for the $G$-action on $X$. 
In particular, $X/G$ and $Z = (X/H)/G$ have the same
dimension. By Definition~\ref{def1} we only need 
to show that

\smallskip
(i) given a canonical form map $F \colon X \dasharrow X$
there exists a canonical form map 
$F' \colon X/H \dasharrow X/H$ such that $\dim \, F(X)
\ge \dim \, F'(X)$ and conversely, that

\smallskip
(ii) given a canonical form map $F' \colon X/H \dasharrow X/H$,
there exists a canonical form map $F \colon X \dasharrow X$ 
such that $\dim \, F(X) = \dim \, F'(X)$. 

\smallskip
Since $H$ is special and its action on $X$ 
is generically free, we can choose a rational 
section $\alpha \colon X/H \dasharrow X$ for the quotient 
map $\pi \colon X \lra X/H$. 

\smallskip
We now proceed to prove (i). 
Suppose $F$ is induced by a rational 
map $f \colon X \dasharrow G$. After translating $\alpha$ by 
an element of $H$, we may assume that the image 
of $\alpha$ meets the domain of $f$. Now let 
$F' \colon X/H \dasharrow X/H$ be the canonical form map induced by
$f' = f \circ \alpha \colon X/H \dasharrow G$. Then 
the diagram
\begin{equation} \label{e.prop2-4}
\xymatrix{X \ar@{-->}[r]^{F} \ar@{->}[d]^{\pi} & X \ar@{->}[d]^{\pi} \\
X/H \ar@{-->}[r]^{F'} & X/H } \end{equation}
commutes, and (i) follows.

To prove (ii), choose a Zariski open subset 
$U \subset X/H$, such that $\alpha$ is regular in $U$ 
and $\pi \colon \pi^{-1}(U) \lra U$ is an $H$-torsor. 
In particular, there is a morphism
$s \colon \pi^{-1}(U) \lra H$ such that 
$s(x) \cdot x = \alpha(\pi(x))$.  
Suppose the canonical form map
$F' \colon X/H \dasharrow X/H$ is induced by
$f' \colon X/H \dasharrow G$.  After replacing $f'$ 
by $g f'$ for a suitable $g \in G$, we may assume 
without loss of generality (and without changing $\dim \, F'(X/H)$)
that $F'(\overline{x}) \in U$ for
$\overline{x} \in X/H$ in general position.

We will now construct the canonical form 
map $F \colon X \dasharrow X$, whose existence is
asserted by (ii). To motivate this construction, 
we remark that it is easy to define a canonical 
form map $F$ so that the diagram 
\eqref{e.prop2-4} commutes; such a map is induced
by $f = f' \circ \pi \colon X \dasharrow G$.  However, 
this diagram only shows that $\dim \, F(X) \ge \dim \, F'(X/H)$;
equality will not hold in general. On the other hand, we are
free to modify $f(x)$ by multiplying $f'(\pi(x))$ by 
any element of $H$ on the left (this element of $H$ may even 
depend on $x$); the resulting canonical form map $F$ 
will still give rise to a commutative 
diagram~\eqref{e.prop2-4}. With this in mind,
we define $F \colon X \dasharrow X$ as the canonical form 
map induced by $f \colon X \dasharrow G$, where
\[ f(x) = s(f'(\pi(x)) \cdot x) f'(\pi(x)) \]
for $x \in X$ in general position. As we mentioned above,
the diagram~\eqref{e.prop2-4} commutes; 
moreover, by our choice of $s$,
$F(X) \subset \alpha(X/H)$. Hence, $\pi$ restricts 
to an isomorphism between $F(X)$ and 
$\pi(F(X)) = F'(X/H)$. This proves (ii).
\end{proof}

\begin{remark} \label{rem.normal}
Lemma~\ref{lem3} and Proposition~\ref{prop.product2}
may fail if $H$ is not special; see Example~\ref{ex1.hom-space}. 
If $H$ is special, Proposition~\ref{ex1.hom-space} may still fail if 
the $H$-action on $X$ is not generically free; 
see Remark~\ref{rem.product2}.
\end{remark}

\section{A lower bound}
\label{sect.lower}

\begin{defn} \label{def.e}
Let $S$ be an algebraic group and $Y$ be a generically free
$S$-variety. We define $e(Y, S)$ as the smallest integer $e$ with
the following property: given a point $y \in Y$ in general position,
there is an $S$-equivariant rational map $f \colon Y \dasharrow Y$
such that $f(Y)$ contains $y$ and $\dim \, f(Y) \le e + \dim(S)$.
\end{defn}

\begin{remark} \label{rem.e05}
Note that this definition is similar to the definition of
the essential dimension $\ed(Y, S)$ of $Y$; cf. Section~\ref{prel.ed}.
The difference is that $\ed(Y, S)$ is the minimal value
of $\dim \, f(Y) - \dim(S)$,
where $f$ is allowed to range over a wider class of rational
$S$-equivariant maps. In particular, $e(Y, S) \ge \ed(Y, S)$.
Note also that $e(Y, S)$ depends only
on the birational class of $Y$, as an $S$-variety.
\end{remark}

\begin{remark} \label{rem.e-grp}
In the sequel we will be particularly interested in the case
where $Y$ is itself an algebraic group, $S$ is a closed subgroup 
of $Y$, and
the $S$-action on $Y$ is given by translations (say, by right 
translations, to be precise).  In this situation, $e(Y, S)$ is simply 
the minimal possible value of $\dim \, f(Y) - \dim(S)$, where 
$f$ ranges over
all $S$-equivariant rational maps $Y \dasharrow Y$. Indeed, 
after composing $f$ with a suitable left translation $g \colon Y \lra Y$, 
we may assume that $f(Y)$ contains any given $y \in Y$. 
\end{remark}

\begin{lem} \label{lem.e1} Let $Y$ be a generically free $S$-variety.

\smallskip
(a) If $Y$ is split (cf. Section~\ref{prel.split}) then $e(Y, S) = 0$.

\smallskip
(b) Suppose there exists a dominant rational $S$-equivariant
map $\alpha \colon V \dasharrow Y$, where $V$ is a vector space
with a linear $S$-action.  Then $e(Y, S) = \ed(S)$.

\smallskip
(c) If $Y = G$ is a special algebraic group, $S$ is a subgroup of $G$ and
the $S$-action on $Y$ is given by translations then $e(Y, S) = \ed(S)$.
\end{lem}

Note that the condition of part (a) is always satisfied if
$S$ is a special group.

\begin{proof}
(a) If $Y$ is split, it is birationally isomorphic to $S \times Z$,
where $S$ acts by translations on the first factor and trivially
on the second. In fact, we may assume without loss of generality
that $Y = S \times Z$.  Now for any $z_0 \in Z$
consider $f_{z_0} \colon S \times Z \dasharrow S \times Z$,
given by $(s, z) \mapsto (s, z_0)$. As $z_0$ ranges over $Z$,
the images of $f_{z_0}$ cover $Y$. Each of these images
has the same dimension as $S$; this yields $e(Y, S) = 0$.

\smallskip
(b) Let $\beta \colon Y \dasharrow Y_0$ be the dominant $S$-equivariant
rational map from $Y$ to a generically free $S$-variety $Y_0$ of minimal
possible dimension, $\ed(S) + \dim(S)$; cf. Remark~\ref{rem.e-grp}.
Then for any $v \in V$,
there is a rational
$G$-equivariant map $\gamma \colon Y_0 \dasharrow V$ such that $v$ lies in
the image of $\gamma$; see~\cite[Proposition 7.1]{reichstein}.
Taking $f = \alpha \circ \gamma \circ \beta \colon Y \dasharrow Y$
in Definition~\ref{def.e} and varying $v$ over $V$,
we see that $e(Y, S) \le \dim(Y_0) - \dim(S) = \ed(S)$.
The opposite inequality was noted in Remark~\ref{rem.e05}.

\smallskip
(c) Let $V$ be a generically free linear representation 
of $G$ (and thus of $S$).  Since $G$ is special, 
$V$ is split; cf. Section~\ref{prel.split}.  Consequently, 
there is a dominant rational map $V \dasharrow G$
of $G$-varieties (and hence, of $S$-varieties). The desired 
conclusion now follows from part (b).
\end{proof}

\begin{prop} \label{prop.lower1} Let $G$ be a connected group and
$X$ be an irreducible $G$-variety with a stabilizer $S$
in general position. Then

\smallskip
(a) $\cd(X, G) \ge e(G, S)$, where $S$ acts on $G$ by translations.

\smallskip
\noindent
In particular,

\smallskip
(b) $\cd(X, G) \ge \ed(G, S)$,
and

\smallskip
(c) if $G$ is special then $\cd(X, G) \ge \ed(S)$.
\end{prop}

\begin{proof} (b) and (c) follow from (a) by Remark~\ref{rem.e05}
and Lemma~\ref{lem.e1}(c) respectively.

To prove part (a), choose a canonical form map
$F \colon X \dasharrow X$ 
such that $\dim \, F(G \cdot x) = \cd(X, G)$
for $x$ in general position; cf. Lemma~\ref{lem2}.
Suppose that $F$ is induced by
a rational map $f \colon X \dasharrow G$, as in Definition~\ref{def.fandF},
and consider the commutative diagram  
$$ \xymatrix{G \ar@{-->}[r]^{F'} \ar@{->}[d]^{\phi} & G \ar@{->}[d]^{\phi} \\
G \cdot x \ar@{-->}[r]^{F} & G \cdot x }$$
of rational maps, where $\phi \colon G \lra G\cdot x$ 
is the orbit map, $\phi(g) = g \cdot x$,
and $F'(g) = f(g \cdot x) \cdot x$.

Now set $S = \Stab_G(x)$ and observe that 
$F'(g s^{-1}) = F'(g) s^{-1}$ for every $s \in S$.
In view of Remark~\ref{rem.e-grp}, this implies
\[ \dim \, F'(G) - \dim \, S \ge e(G, S) \, . \]
On the other hand, $F'(G)$ is an $S$-invariant subvariety 
of $G$ and (because the above diagram in commutative)
$F(G \cdot x) = \phi (F'(G))$. Finally,
since for any $g \in G$, 
$\phi(g s^{-1}) = \phi(g)$ if and only if
$s \in S$, we see that the fibers of the map
$\phi_{| F'(X)} \colon F'(X) \lra F(G \cdot x)$ 
are precisely the $S$-orbits in $F'(X)$.
Consequently, 
\[ \cd(X, G) = \dim \, F(G \cdot x) = 
\dim \, F'(G) - \dim \, S \ge e(G, S) \, , \]
as claimed.
\end{proof}

\begin{remark} \label{rem3.2}
Proposition~\ref{prop.lower1} assumes
that the $G$-action on $X$ has a stabilizer 
in general position, i.e., there exists 
a subgroup $S \subset G$
such that $\Stab(x)$ is conjugate 
to $S$ for $x \in X$ in general position.
This condition is satisfied by many but 
not all group actions; see~\cite[Section 7]{pv}. 
For an arbitrary $G$-action on $X$, 
our proof of part (a) shows that if
$e(G, \Stab(x)) \ge d$ for $x$ in
a Zariski dense open subset of $X$ then $\cd(X, G) \ge d$.

Note also that if $L_x$ is a Levi subgroup of $\Stab(x)$
then by a theorem of Richardson (see~\cite[Theorem 9.3.1]{richardson}
or \cite[Theorem 7.1]{pv}), there exists a non-empty 
Zariski open subset $U \subset X$ such that $L_x$ and $L_y$ 
are conjugate in $G$. Since $\ed(L_x) = \ed(\Stab(x))$
(this is an immediate consequence of \cite[Lemma 1.13]{sansuc}; 
for a direct geometric proof, see~\cite{kordonskii}), 
$\ed(\Stab(x))$ assumes the same value for every $x \in U$.
In particular, Proposition~\ref{prop.lower1}(c)
remains valid for an arbitrary $G$-action, provided that 
we replace the inequality $\cd(X, G) \ge \ed(S)$ by
$\cd(X, G) \ge \ed(\Stab(x))$ for $x \in U$.
\end{remark}

\begin{cor} \label{cor4.3} Let $G$ be a connected group, $S$ be
a closed subgroup, and $X = G/S$ be a homogeneous space.
Then

\smallskip
(a) $\cd(X, G) = e(G, S)$, where $S$ acts on $G$ by translations.

\smallskip
(b) If $G$ is special then $\cd(X, G) = \ed(S)$.
\end{cor}

\begin{proof} Part (b) follows from part (a) and Lemma~\ref{lem.e1}(c).

To prove (a), note that by Proposition~\ref{prop.lower1},
we only need to show that
$\cd(X, G) \le e(G, S)$, i.e., to construct a canonical form map
$F \colon X \dasharrow X$ such that
\begin{equation} \label{e.F}
\dim \, F(G \cdot x) = e(G, S)
\end{equation}
for $x$ in general position.  We will define $F$
by reversing the construction in Proposition~\ref{prop.lower1}.
Let $F' \colon G \dasharrow G$ be an $S$-equivariant rational map
(with respect to the right translation action of $S$ on $G$), such that
$\dim \, F'(G)$ assumes its minimal possible value, 
$e(G, S) + \dim(S)$; cf. Remark~\ref{rem.e-grp}.
Then $f' \colon G \dasharrow G$ given by $f'(g) = F'(g)g^{-1}$
is $S$-invariant (with respect to the right translation action 
of $S$ on $G$).  Hence, $f'$ descends to $f \colon G/S \dasharrow G$.
Thus we have a commutative diagram
\[   \begin{array}{ccc}
  G & \stackrel{F'}{\dasharrow} & G  \\
  \downarrow &  & \downarrow  \\
  G/S & \stackrel{F}{\dasharrow} & G/S
\end{array} \]
where $F(x) = f(x) \cdot x$. Here $F$ is, by construction, a canonical
form map, and
\[ \dim \, F(G/S) = \dim \, F'(G) - \dim \, S = e(G, S) \, , \]
as desired.
\end{proof}

\begin{example} \label{ex1.hom-space} 
Let $G$ be a special group and $H$ be a non-special
closed normal subgroup of $G$. (For example, 
$G = \SL_n$ and $H = \mu_n$ is the center of $G$.)
 
\smallskip
(i) Let $X = G/H$. Then by Corollary~\ref{cor4.3}(b), 
$\cd(X, G) = \ed(H)$, which is $\ge 1$;
cf. \cite[Proposition 5.3]{reichstein}.
On the other hand, $\cd(X, G/H) = 0$; 
(cf. Lemma~\ref{lem2.2}(a)). This shows that
the equality $\cd(X, G) = \cd(X, G/H)$
in Lemma~\ref{lem3}(b) may fail if $H$ is not special.
 
\smallskip
(ii) Now let $X = G$ (viewed as a $G$-variety with the translation action).
Then $\cd(X, G) = 0$ (cf. Lemma~\ref{lem2.2}(a)) 
but,  $\cd(X/H, G) = \ed(H)$.
This shows that the equality $\cd(X, G) = \cd(X/H, G)$ in
Proposition~\ref{prop.product2} may also fail if $H$ is not special.
\end{example}

\section{A comparison lemma}
\label{sect.rel-sect}

\begin{lem} \label{lem5.2}
Let $\alpha \colon X \dasharrow Y$ be a dominant rational
map of irreducible $G$-varieties. Suppose 
$\dim(G \cdot x) = d$ and $\dim(G \cdot y) = e$
for $x \in X$ and $y \in Y$ in general position.
Then
$\cd(X, G) \le \cd(Y, G) + d - e$.
\end{lem}

\begin{proof} Let $F \colon Y \dasharrow Y$ be a canonical form map
such that  $\dim \, F(G \cdot y) = \cd(Y, G)$ for $Y$ in general position. 
Suppose $F$ is induced by $f \colon Y \dasharrow G$, i.e., 
$F(y) = f(y) \cdot y$, as in Definition~\ref{def.fandF}.

Now consider $f' = f \circ \alpha \colon X \dasharrow G$ and the
induced canonical form map $F' \colon X \dasharrow X$ given by
$F'(x) = f'(x) \cdot x$.  The relationship between $F$
and $F'$ is illustrated by the following commutative diagram,
where $x$ is a point in general position in $X$ 
and $y = \alpha(x) \in Y$.
\[   \begin{array}{ccc}
  G \cdot x & \stackrel{F'}{\dasharrow} & F'(G \cdot x)  \\
\alpha \; \downarrow \; \; \; &  & \; \; \; \downarrow  \; \alpha \\
  G \cdot y & \stackrel{F}{\dasharrow} & F(G \cdot y) \, .
\end{array} \]
Each fiber of $\alpha \colon G \cdot x \lra G \cdot y$ 
has dimension $d - e$. Hence, each fiber of the right
vertical map $\alpha_{| F'(G \cdot x)}$ has dimension $\le d - e$.  Applying
the fiber dimension theorem to this map, we obtain
\[
\dim \, F'(G \cdot x) \le \dim  \, F(G \cdot y) + d - e = \cd(Y, G) + d - e \, ,
\]
and the proposition follows; cf. Lemma~\ref{lem2}.
\end{proof}

Let $X$ be a $G$-variety and $H$ be a closed subgroup of $G$.
Recall that an $H$-invariant (not necessarily irreducible)
subvariety $Y \subset X$ is called a
$(G, H)$-section if (i) $G \cdot Y$ is dense in $X$ and (ii) for $y \in Y$
in general position, $g \cdot y \in Y \Leftrightarrow g \in H$. Note that
in some papers a $(G, H)$-section is called a {\em relative section}
(cf. \cite[Section 2.8]{pv}) or a {\em standard relative section
with normalizer} $H$ (cf~\cite[(1.7.6)]{popov}).

\begin{cor} \label{cor5.3} Let $X$ be an irreducible $G$-variety.

\smallskip
(a) If $X$ has a $(G, H)$-section then
$\cd(X, G) \le e(G, H) + d - \dim(G) + \dim(H)$, where
$d = \dim(G \cdot x)$ for $x \in X$ in general position.

\smallskip
(b) If $X$ has a stabilizer $S$ in general position then
\[ e(G, S) \le \cd(X, G) \le e(G, N) - \dim(S) + \dim(N) \, , \]
where $N$ is the normalizer of $S$ in $G$.
\end{cor}

\begin{proof} (a) The existence of a $(G, H)$-section is equivalent to the
existence of a $G$-equivariant rational map $X \dasharrow G/H$;
see~\cite[Theorem 1.7.5]{popov}. Thus by Lemma~\ref{lem5.2},
$\cd(X, G) \le \cd(G/H, G) - d + \dim(G/H)$. By
Corollary~\ref{cor4.3}(a) $\cd(G/H, G) = e(G, H)$, and part (a) follows.

\smallskip
(b) The inequality $e(G, S) \le \cd(X, G)$ follows
from Proposition~\ref{prop.lower1}(a).
To prove the inequality
\begin{equation} \label{e.cor5.3}
\cd(X, G) \le e(G, N) - \dim(S) + \dim(N) \, ,
\end{equation}
note that by \cite[(1.7.8)]{popov}, $X$ has
a $(G, N)$-section. Substituting $H = N$ and $d = \dim(G) - \dim(S)$
into the inequality of part (a), we obtain~\eqref{e.cor5.3}.
\end{proof}

\section{The canonical dimension of a group}
\label{sect.defn-grp}

In this section we will define the canonical dimension of an algebraic
group $G$. We begin with a simple lemma.

\begin{lem} \label{lem5.1} Let $X$ be an irreducible
$G$-variety, and let $Z$ be
an irreducible variety with trivial action of $G$. Then
$\cd(X \times Z, G) = \cd(X, G)$.
\end{lem}

\begin{proof} The inequality $\cd(X \times Z, G) \le \cd(X, G)$ 
follows from Lemma~\ref{lem5.2}, applied to the projection
map $\alpha \colon X \times Z \lra X$. To prove the opposite 
inequality, let $c = \cd(X \times Z, G)$ and
choose a canonical form map $F \colon X \times Z \dasharrow X \times Z$ 
such that $\dim \, F(G \cdot (x, z)) = c$. Suppose $F$ is induced by
$f \colon X \times Z \dasharrow G$, i.e., 
$F(x, z) = f(x, z) \cdot (x, z)$, as in Definition~\ref{def.fandF}.
It is now easy to see that for $z_0 \in Z$ in general
position, the map $f_{z_0} \colon X \dasharrow G$ given by
$f_{z_0}(x) = f(x, z_0)$ gives rise to a canonical form map
$F_{z_0} \colon X \dasharrow X$ such that $\dim \, F_0(G \cdot x) = c$.
In other words, $\cd(X, G) \le c$, as claimed.
\end{proof}

\begin{prop} \label{prop5.3}
Let $V$ be a generically free linear representation of $G$.

\smallskip
(a) If $X$ is an irreducible generically free $G$-variety then
$\cd(X, G) \le \cd(V, G)$.

\smallskip
(b) If $W$ is another generically free $G$-representation, then
$\cd(V, G) = \cd(W, G)$.
\end{prop}

\begin{proof} (a) By \cite[Corollary 2.20]{reichstein},
there is a dominant rational map
$\alpha \colon X \times \bbA^d \dasharrow V$ of $G$-varieties,
where $d = \dim(V)$, and $G$ acts trivially on $\bbA^d$. Now
\[ \cd(X, G) \stackrel{\text{\tiny by Lemma~\ref{lem5.1}}}{=}
\cd(X \times \bbA^d, G)
\stackrel{\text{\tiny by Lemma~\ref{lem5.2}}}{\le}
\cd(V, G) \, , \]
as claimed.

(b) $\cd(W, G) \le \cd(V, G)$ by part (a). To prove the opposite
inequality, interchange the roles of $V$ and $W$.
\end{proof}

\begin{defn} \label{def.cd}
We define the canonical dimension $\cd(G)$ of an algebraic
group $G$ to be $\cd(V, G)$, where $V$ is a generically
free linear representation of $G$.  By Proposition~\ref{prop5.3}(a)
this number is independent of the choice of $V$. Moreover, by
Proposition~\ref{prop5.3}(b)
$\cd(G) = \max \{ \cd(X, G) \}$, as $X$ ranges over
all irreducible generically free $G$-varieties.
\end{defn}

\begin{cor} \label{cor6.05}
Suppose $W$ is a linear representation of $G$ such that
$\Stab_G(w)$ is finite for $w \in W$ in general position. Then
$\cd(G) \le \cd(W, G)$.
\end{cor}

\begin{proof} Let $V$ be a generically free linear representation of
$G$. Then so is $X = V \times W$. The desired inequality is now
a consequence of Lemma~\ref{lem5.2}, applied to the projection
map $\alpha \colon V \times W \lra W$.
\end{proof}

\begin{lem} \label{lem6.1} 
(a) $\cd(G) \le \cd(H) + \dim(G) - \dim(H)$, for any closed
subgroup $H \subset G$.

\smallskip
(b) $\cd(G) = \cd(G^0)$.

\smallskip
(c) $\cd(G) = 0$ if and only if $G^0$ is special.

\smallskip
(d) $\cd(G_1 \times G_2) \le \cd(G_1) + \cd(G_2)$.
\end{lem}

\begin{proof} 
(a) Follows from Lemma~\ref{lem2.1}, with
$X = V$ = generically free linear representation of $G$.

\smallskip
(b) Immediate from Lemma~\ref{lem2.15}.

\smallskip
(c) By part (b), we may assume $G = G^0$ is connected.
The desired conclusion now follows from Lemma~\ref{lem2.2}.

\smallskip
(d) Follows from Lemma~\ref{lem.product}, by taking $X_i = V_i$
to be a generically free representation of $G_i$ for $i = 1, 2$.
\end{proof}

\begin{example} \label{ex6.2}
Consider the subgroup
\[ H = \{ \begin{pmatrix} A & \begin{matrix} b_1 \\ \vdots\\ b_{n-1}
                              \end{matrix} \\
                0 \dots 0  &   1      \end{pmatrix}
\; | \;  A \in \GL_{n-1} \, , \, b_1, \dots, b_{n-1} \in k \} \, , \]
of $G = \PGLn$. The Levi subgroup of $H$ 
is special (it isomorphic to $\GL_{n-1}$); hence,
$H$ itself is special. (This follows from the theorem of Grothendieck
stated in Section~\ref{prel.special} or alternatively, 
from~\cite[Lemma 1.13]{sansuc}.) Since $\dim(H) = n^2-n$,
Lemma~\ref{lem6.1}(a) yields $\cd(\PGL_n) \le n-1$. In particular,
$\cd(\PGL_2) = 1$. (Note that $\cd(\PGL_2) \ge 1$ by Lemma~\ref{lem6.1}(c)).

Alternatively, we can deduce the inequality $\cd(\PGL_n) \le n-1$
by applying Lemma~\ref{lem5.2} to
the projection map $\Mat_n \times \Mat_n \lra \Mat_n$
to the first factor, where $\PGLn$ acts on $\Mat_n$ by conjugation.
The $\PGLn$-action on $\Mat_n \times \Mat_n$ is generically free; hence,
$\cd(\PGLn) = \cd(\Mat_n \times \Mat_n, \PGLn)$. On the other hand,
$\cd(\Mat_n, \PGLn) = 0$; see Example~\ref{ex1}.
Now Lemma~\ref{lem5.2} tells us that
\[ \cd(\PGLn) = \cd(\Mat_n \times \Mat_n, \PGLn) \le
\cd(\Mat_n, \PGLn) + n-1 - 0 = n - 1 \, . \]
For a third proof of this inequality, see Example~\ref{ex6.2a}.
\end{example}

\section{Splitting fields}
\label{sect.spl-fld}

Throughout this section we will assume that $G/k$ is a connected
linear algebraic group. Unless otherwise specified, the fields
$E$, $K$, $L$, etc., are assumed to be finitely 
generated extensions of the base field $k$.

Let $X$ be a generically free irreducible $G$-variety,
$E = k(X)^G = k(X/G)$,
$\pi \colon X \dasharrow X/G$ be the rational quotient map
and $F \colon X \dasharrow X$ be a canonical form map.
Recall that $F$ commutes with $\pi$, so that
$F(X)$ may be viewed as an algebraic variety over $E$.

\begin{lem} \label{lem.split} 
Let $X$ be a generically free $G$-variety such that
$G$-orbits in $X$ are separated by regular invariants and
let $F \colon X \dasharrow X$ be a canonical form map.
Suppose $\alpha \in H^1(E, G)$ is the class represented by $X$.
Then for any field extension $K/E$ 
the following conditions are equivalent: 

\smallskip
(a) $\alpha_K = 1$,

\smallskip
(b) $X$ is rational over $K$,

\smallskip
(c) $F(X)$ is unirational over $K$,

\smallskip
(d) $K$-points are dense in $F(X)$,

\smallskip
(e) $F(X)$ has a $K$-point.
\end{lem}

\begin{proof} We begin by proving the lemma in the case where $K = E$. 

(a) $\Rightarrow$ (b):
If $\alpha = 1$ then $X$ is birationally isomorphic
to $X/G \times G$ (over $X/G$). 
Now recall that the underlying variety of a connected 
algebraic group $G$ is rational over $k$.
Hence, $X/G \times G$ is rational over $X/G$, i.e., $X$ is rational 
over $E$.

(b) $\Rightarrow$ (c): The rational map $F \colon X \dasharrow F(X)$ is,
by definition, dominant. If $X$ is rational, this makes $F(X)$ unirational.

(c) $\Rightarrow$ (d) and (d) $\Rightarrow$ (e) are obvious.

(e) $\Rightarrow$ (a): 
An $E$-point in $F(X)$ is a rational section
$s \colon X/G \dasharrow F(X) \subset X$ for $\pi$.
The existence of such a section implies that $X$ is split, 
and hence, so is $\alpha$; see Section~\ref{prel.split}.

\smallskip
To prove the general case, note that since the $G$-orbits in $X$
are separated by regular invariants, we can choose 
a regular model of the rational quotient variety
$X/G$, so that the rational quotient map $\pi \colon X \lra X/G$
is regular and its fibers are exactly the $G$-orbits in $X$.
After making $X/G$ smaller if necessary, we may also assume that
our field extension $K/E$ is represented by 
a surjective morphism $Y \lra X/G$ of algebraic
varieties. Then $\alpha_K$ is represented
by the $G$-variety $X_K = X \times_{X/G} Y$.

We claim that the morphism $\pi_K \colon X_K \lra Y$ (projection 
to the second component) separates the $G$-orbits
in $X_K$. Indeed, if for some $x_1, x_2 \in X$,
$z_1 = (x_1, y)$ and $z_2 = (x_2, y) \in X_K$ have the same 
second component then $\pi(x_1) = \pi(x_2)$. We conclude that
$x_1$ and $x_2$ are in the same $G$-orbit in $X$ and consequently,
$z_1$ and $z_2$ are in the same $G$-orbit in $X_K$, as claimed.
This shows that $\pi_K$ is a rational quotient map for the $G$-action
on $X_K$; cf. \cite[Lemma 2.1]{pv}.

We now define a rational map
$F_K \colon X_K \dasharrow X_K$ by $F_K (x, y) = (F(x), y)$.
Since $\pi_K \circ F_K = \pi_K$, $F_K$ is a canonical form map;
see Remark~\ref{rem1}. Moreover, $F_K(X_K) = F(X) \times_{X/G} Y$.
Replacing $X$ by $X_K$ and $F$ by $F_K$, we reduce 
the lemma to the case we settled at the beginning of the proof
(where $K = E$). 
\end{proof}

Let $\alpha \in H^1(E, G)$.  As usual, we will
call a field extension $K/E$ a {\em splitting field} for $\alpha$
if the image $\alpha_K$ of $\alpha$
under the natural map $H^1(E, G) \lra H^1(K, G)$ is split.
If $\alpha$ is represented by a generically free $G$-variety $X$,
with $k(X)^G = E$ then we will also sometimes say that $K$ is a splitting
field for $X$.

\begin{prop} \label{prop.spl-deg}
Suppose $\alpha \in H^1(E, G)$.

\smallskip
(a) If $\cd(\alpha) = 1$ then there exist $0 \ne a, b \in E$ such 
that a field extension $K/E$ splits $\alpha$ if and only 
if the quadratic form
$q(x, y, z) = x^2 + a y^2 + bz^2$ is isotropic over $K$.
In particular, $\alpha$ has a splitting field
$K/E$ of degree $2$.

\smallskip
(b) If $\cd(\alpha) = 2$ then $\alpha$ has a splitting field
$K/E$ of degree $2$, $3$, $4$, or $6$.
\end{prop}

Note that if $K/E$ is a splitting field for $\alpha$
then $[K:E] = 1$ is impossible in either part.  Indeed, 
otherwise $\alpha$ itself is split, and $\cd(\alpha) = 0$ by
see Lemma~\ref{lem2.2}(a).

\begin{proof}
Choose a canonical form map $F \colon X \dasharrow X$, such that 
$\dim \, F(X) - \dim(X/G) = \cd(\alpha)$. By Lemma~\ref{lem.split},  
$F(X)$ is unirational over every splitting field $K$ of $\alpha$; 
in particular, it is unirational over the algebraic
closure $\overline{E}$ of $E$. 

\smallskip
(a) Here $F(X)$ is a curve over $E$, and L\"uroth's theorem 
tells us that $F(X)$ is rational over $\overline{E}$.  
It is well known that any such curve is birationally isomorphic 
to a conic $Z$ in $\bbP_E^2$ (see, e.g.,
\cite[Proposition 1.1.1]{mt}) and that $K$-points are dense in $Z$
if and only if $Z(K) \ne \emptyset$ (see,
e.g.,~\cite[Theorem 1.2.1]{mt}).  Writing the equation of
$Z \subset \bbP_E^2$ in the form $x^2 + ay^2 + bz^2 = 0$,
we deduce the first assertion of part (a). The second assertion 
is an immediate consequence of the first; for example, 
$K = E(\sqrt{-a})$ is a splitting field for $\alpha$.

\smallskip
(b) Here $F(X)$ is a surface over $E$, which becomes unirational
over the algebraic closure $\overline{E}$.
By a theorem of Castelnuovo, $F(X)$ is, in fact, rational 
over $\overline{E}$. Let $Z$ be a complete smooth minimal surface,
defined over $E$, which is birationally isomorphic to $F(X)$ 
via $\phi \colon F(X) \stackrel{\sim}{\dasharrow} Z$ and
let $U \subset Z$ be an open subset such that
$\phi$ is an isomorphism over $U$.  Part (b) now follows 
from Lemma~\ref{lem.split} and Lemma~\ref{lem.isk} below.  
\end{proof}

\begin{lem} \label{lem.isk}
Let $E$ be a field of characteristic zero,
$Z$ be a complete minimal surface defined over $E$ 
and rational over $\overline{E}$,
and let $U$ be a dense open subset of $Z$ (defined over $E$).
Then $U$ contains a $K$-point for some field extension
$K/E$ of degree $1$, $2$, $3$, $4$ or $6$.
\end{lem}

\begin{proof} 
By the Enriques-Manin-Iskovskih classification 
$Z$ is a conic bundle or a del Pezzo surface; 
see~\cite[Theorem 1]{iskovskih} or \cite[Theorem 3.1.1]{mt}.  
Note that $Z = \bbP^2$, listed as a separate 
case in~\cite[Theorem 1]{iskovskih}, is, in fact, 
a del Pezzo surface. (We remark however, that 
the lemma is obvious in this case, since $E$-points 
are dense in $\bbP_E^2$.)

If $f \colon Z \lra C$ is a conic bundle over a rational
curve $C$, then after replacing $E$ by a quadratic
extension $E'$, we may assume that $C_{E'} \simeq \bbP_{E'}^1$.
For every $E'$-point $z \in C$, $f^{-1}(z)$ is
a rational curve over $E'$. Taking $z \in C_{E'}$ so that
$f^{-1}(z) \cap U \ne \emptyset$, we can choose an extension
$K/E'$ of degree $1$ or $2$ so that $f^{-1}(z)_K \simeq \bbP_K^1$.
Now $[K:E] = 1$, $2$ or $4$, and $K$-points are dense in $f^{-1}(z)$,
so that one of them will lie in $U$.

From now on we may assume that $Z$ is a del Pezzo surface.  Recall that
the anticanonical divisor $-\Omega_Z$ on a del Pezzo surface
is ample, and the degree $d = \Omega_Z \cdot \Omega_Z$ can range
from $1$ to $9$.

If $d = 1$ the linear system $| -2 \Omega_Z|$ defines a (ramified) double
cover $f \colon Z \lra Q$, where $Q$ is a quadric cone in $\bbP_E^3$;
see~\cite[p. 30]{iskovskih}. Then $Q_{E'} \simeq \bbP_{E'}^2$ for
some extension $E'/E$ of degree $1$ or $2$. Now choose an $E'$-point
$x \in f(U) \subset Q$ and split $f^{-1}(x)$ over a field
extension $K/E'$ of degree $1$ or $2$. Then $[K:E] = 1$, $2$ or $4$
and $U$ contains a $K$-point.

If $d = 2$ then the linear system $|- \Omega_Z|$ defines a (ramified)
double cover $Z \lra \bbP_E^2$ (see~\cite[p. 30]{iskovskih}),
and points of degree $2$ are dense in $Z$.

If $3 \le d \le 9$ then it is enough to show that $Z(K) \ne \emptyset$
for some field extension $K/E$ of degree $1$, $2$, $3$, $4$ or $6$.
Indeed, if $Z(K) \ne \emptyset$ then $Z_K$ is unirational over $K$
(see~\cite[Theorem 3.5.1]{mt}) and thus $K$-points are dense in $Z$.
Note also that for $3 \le d \le 9$,
$Z$ is isomorphic to a surface in $\bbP^d$ of
degree $d$. Intersecting this surface by two hyperplanes
in general position, we see that $Z$ has a point of degree
dividing $d$. This proves the lemma for $d = 3$, $4$, and $6$.

For $d = 5$ and $7$, $Z$ always has an $E$-point
(see~\cite[Theorem 7.1.1]{mt}), so the lemma holds trivially
in these cases. For $d = 8$, $Z$ has a point of degree dividing $4$ and
for $d = 9$, $Z$ has a point of degree dividing $3$ (see \cite[p. 80]{mt}).
The proof of the lemma is now complete.
\end{proof}

\begin{example} \label{ex.pgl}
Suppose $\alpha \in H^1(K, \PGL_n)$ is represented by a
central simple algebra of index $d$. Then the degree
of every splitting field for $\alpha$ is divisible
by $d$ (cf. e.g, \cite[Theorem 7.2.3]{rowen}); hence,
$\cd(\alpha) \ge \begin{cases} \text{$2$, if $d \ge 3$,} \\
\text{$3$, if $d \ne 1, 2, 3, 4$ or $6$.} \end{cases}$.
In particular,
$\cd(\PGLn) \ge \begin{cases} \text{$2$, if $n \ge 3$,} \\
\text{$3$, if $n \ne 1, 2, 3, 4$ or $6$.} \end{cases}$
For sharper results on $\cd(\PGLn)$, see Section~\ref{sect.a_{n-1}}.  
\end{example}

\begin{example} \label{ex.exceptional1}
Let $V$ be a generically free linear representation of
$G = F_4$,  $E_6$ or $E_7$ (adjoint or simply connected),
$K = k(V)^G$ and $\alpha \in H^1(K, G)$ be the class represented
by the $G$-variety $V$.  Then the degree of any splitting field
$L/K$ for $\alpha \in H^1(k(V)^G, G)$ is divisible
by $6$;~\cite[p. 223]{ry2}.  We conclude that
$\cd(G) \ge 2$ for these groups.
\end{example}

\begin{example} \label{ex.exceptional2}
$\cd(G) \ge 3$, if $G = E_8$ or adjoint $E_7$;
see~\cite[Corollaries 5.5 and 5.8]{ry2}.
\end{example}

\begin{remark}
Let $G$ be a connected linear algebraic group defined over $k$
and let $H$ be a finite abelian $p$-subgroup of $G$, where $p$
is a prime integer. Recall that the {\em depth} of $H$
is the smallest value of $i$ such that $[H : H \cap T] = p^i$,
as $T$ ranges over all maximal tori of $G$;
see~\cite[Definition 4.5]{ry2}. Note that $H$ has depth $0$
if and only if it lies in a torus of $G$. A prime $p$ is called
a {\em torsion prime} for $G$ if and only if $G$ has a finite abelian
$p$-subgroup of depth $\ge 1$. (This is one of many equivalent
definitions of torsion primes; see~\cite[Theorem~2.28]{steinberg}.)
The inequalities of
Examples~\ref{ex.pgl} - \ref{ex.exceptional2} may be viewed as
special cases of the following assertion:

{\em
Suppose a connected linear algebraic group $G$ has 
a $p$-subgroup $H$ of depth $d$.

\smallskip
(a) If $\cd(G) \le 1$ then $p^d = 1$ or $2$.

\smallskip
(b) If $\cd(G) \le 2$ then $p^d = 1$, $2$, $3$ or $4$.}

\smallskip
\noindent
The proof is immediate from Proposition~\ref{prop.spl-deg} and
\cite[Theorem 4.7]{ry2} (where we take
$X$ to be a generically free linear representation of $G$).
\end{remark}

\section{Generic splitting fields}
\label{sect.gen-spl}

\begin{defn} \label{def.splitting1}
Let $K/E$ be a (finitely generated) field extension. 
A (finitely generated) field extension $L/E$ is said to be
 
\smallskip
(a) a {\em specialization} of $K/E$ if
there is a place $\phi \colon K \lra L \cup \{ \infty \}$, 
defined over $E$,
 
\smallskip
(b) a {\em rational specialization} of $K/E$ if there
is an embedding $K \hookrightarrow L(t_1, \dots, t_r)$,
over $E$, for some $r \ge 0$.

\smallskip
In geometric language, (a) and (b) can be restated as follows.
Suppose the field extensions $K/E$ and $L/E$ are induced by 
dominant rational maps $V \dasharrow Z$ and $W \dasharrow Z$
of irreducible algebraic $k$-varieties, respectively. Then

\smallskip
(a$'$) $W$ is a specialization of $V$ if there 
is a rational map $W \dasharrow V$, such that the diagram
\[ \xymatrix{W \ar@{-->}[dr] \ar@{-->}[r]
 & V \ar@{-->}[d] \\
& Z  } \]
commutes.

\smallskip
(b$'$) $W$ is a rational specialization of $V$ if for some $r \ge 0$
there is a dominant map $W \times \bbA^r \dasharrow V$ such that
the diagram
\[ \xymatrix{W \times \bbA^r \ar@{-->}[dr] \ar@{-->>}[r]
 & V \ar@{-->}[d] \\ & Z  } \]
commutes.
\end{defn}

\begin{remark} \label{rem.splitting1.1}
In the definition of rational specialization
we may assume without loss of generality that
$r = \max \, \{ 0, \trdeg_E(L) - \trdeg_E(K) \}$;
see~\cite[Lemma 1]{roquette2}.
\end{remark}

\begin{defn} \label{def.splitting1.2} Let $\alpha \in H^1(E, G)$.
A splitting field $K/E$ for $\alpha$ is called {\em generic}
(respectively, {\em very generic}) if every splitting field $L/E$ 
for $\alpha$ is a specialization (respectively, a rational 
specialization) of $K/E$.
\end{defn}

\begin{remark} \label{rem.splitting1.1a}
It is easy to see that a rational specialization is a specialization
(cf.~\cite[Lemma 11.1]{saltman}). Consequently, a very generic 
splitting field is generic.
\end{remark}

\begin{remark} \label{rem.kr}
The generic splitting field of $\alpha$
in Definition~\ref{def.splitting1.2} 
is the same as the generic splitting field for the twisted
group $_{\alpha}G$ defined by Kersten and Rehmann~\cite{kr}.
\end{remark}

\begin{lem} \label{lem.splitting2}
Let $G$ be a connected algebraic group, $X$ be an 
irreducible generically free $G$-variety,
$E = k(X)^G = k(X/G)$ and $F \colon X \dasharrow X$
be a canonical form map.  Then $k(F(X))/E$ is a very generic 
splitting field for the class $\alpha \in H^1(E, G)$ 
represented by $X$. 
\end{lem}

\begin{proof} After replacing $X$ by a $G$-invariant open subset,
we may assume that $G$-orbits in $X$ are separated by 
regular invariants. The generic point of $F(X)$ is 
a $k(F(X))$-point; hence, by Lemma~\ref{lem.split}, $F(X)/E$ 
is a splitting field for $\alpha$.

It remains to show that 
every splitting field $L/E$ for $\alpha$ is a rational specialization
of $k(F(X))/E$.  After replacing $X$ by a smaller $G$-invariant 
dense open subset, we may assume that
$L/E$ is induced by a surjective morphism
$Y \lra X/G$ of algebraic varieties.  Then $\alpha_L = 1_L$ 
is represented by 
the generically free $G$-variety $X_L = X \times_{X/G} Y$. 
Since $X_L$ is split, it is rational over $L$; see Lemma~\ref{lem.split}.
The morphisms 
\[ X_L \stackrel{\operatorname{pr}_1}{\lra} X \stackrel{F}{\lra} F(X) 
\stackrel{\pi}{\lra} X/G \]
now tell us that $k(F(X)) \hookrightarrow k(X_L) = L(t_1, \dots, t_r)$, 
over $E$, where $t_1, \dots, t_r$ are independent variables 
and $r = \dim(G)$. (Here $\operatorname{pr}_1$ is the projection
$X_L = X \times_{X/G} Y \lra X$ to the first factor.) 
This shows that $L/E$ is a rational specialization of $k(F(X))/E$, 
as claimed.
\end{proof}

\begin{prop} \label{prop.splitting4}
Let $E/k$ be a finitely generated field extension. Then
for every $\alpha \in H^1(E, G)$,
\begin{eqnarray*} \text{$\cd(\alpha) = \min \, \{ \trdeg_E(K) \; | \; K/E$
is a generic splitting field for $\alpha \, \}$}  \\
 = \text{$\min \, \{ \trdeg_E(L) \; | \; L/E$
is a very generic splitting field for $\alpha \, \}$.}
\end{eqnarray*}
\end{prop}

\begin{proof} Let $X$ be a generically free 
$G$-variety representing $\alpha$;
in particular, $E = k(X)^G = k(X/G)$.
Since $\cd(X, G)$ is, by definition, the minimal value of
$\dim \, F(X) - \dim(X/G) = \trdeg_E \, k(F(X))$,
as $F$ ranges over all canonical form maps
$X \dasharrow X$, Lemma~\ref{lem.splitting2} tells us that
\[ \text{$\cd(\alpha) = \cd(X, G) \ge \min \{ \trdeg_E(L) \; | \; L/E$
is a very splitting field for $\alpha \, \}$.} \]

Now let $K/E$ be a generic splitting field for $\alpha$.
It remains to show that
\begin{equation} \label{e.splitting5}
\cd(X, G) \le \trdeg_E(K)
\end{equation}
Choose a variety $Y$ whose function field $k(Y)$ is $K$;
the inclusion $E \subset K$ then gives rise to a rational map
$Y \dasharrow X/G$.
By Lemma~\ref{lem.splitting2} (with $F = id$), $k(X)/E$ is a very generic
(and hence, a generic; cf. Remark~\ref{rem.splitting1.1a}) 
splitting field of $\alpha$. Since $k(X)/E$
and $K/E$ are both generic splitting fields, 
each is a specialization of the other. 
Geometrically, this means that there exist rational
maps $f_1 \colon X \dasharrow Y$ and $f_2 \colon Y \dasharrow X$
such that the diagram
$$\xymatrix{X\ar@{-->}[r]^{f_1}\ar@{-->}[dr] & 
Y\ar@{-->}[r]^{f_2}\ar@{-->}[d]&X\ar@{-->}[dl]\\
& X/G & }$$
commutes. After replacing $Y$ by the closure of the graph of $f_2$
in $Y \times X$, and $f_2$ by the projection from this
graph to $X$, we may assume that $f_2$ is a morphism.
Now $F = f_2 \circ f_1$ is a well defined rational 
map $X \dasharrow X$ which commutes
with the rational quotient map $\pi \colon X \dasharrow X/G$.
By Remark~\ref{rem1}, $F$ is a canonical form map. Thus
 \begin{eqnarray*}
 \cd(X, G) \le \dim \, F(X)  - \dim \, X/G \le
 \dim \,  f_2(Y) - \dim \, X/G \le \\
 \dim(Y) - \dim(X/G) = \trdeg_k(K) - \trdeg_k (E) = \trdeg_E(K)
 \end{eqnarray*}
Thus completes the proof of~\eqref{e.splitting5} and thus of
Proposition~\ref{prop.splitting4}.
\end{proof}

\begin{example} \label{ex6.2a}
Let $\alpha \in H^1(E, \PGLn)$ be represented by a central simple
$E$-algebra $A$ and let $K$ be the function field of the
Brauer-Severi variety of $A$. Then $K/E$ is a very generic
splitting field for $\alpha$ (see, e.g.,~\cite[Corollary 13.9]{saltman})
and $\trdeg_E(K) = n-1$. By Proposition~\ref{prop.splitting4},
$\cd(\alpha) \le n-1$. This gives yet another proof of
the inequality $\cd(\PGL_n) \le n-1$ of Example~\ref{ex6.2}.
\end{example}

\section{The canonical dimension of a functor}
\label{sect.cd-functor}

The results of the previous section naturally lead to
the following definitions.  Let $\mathcal{F}$ be
a functor from the category $\Fields$ of finitely generated extensions of
the base field $k$ to the category
$\Sets^*$ of pointed sets.  We will
denote the marked element in $\mathcal{F}(E)$ by $1_E$ (and
sometimes simply by $1$, if the reference to the field $E$ is
clear from the context).  Given a field extension $L/E$,
we will denote the image of $\alpha \in \mathcal{F}(E)$
in $\mathcal{F}(L)$ by $\alpha_L$.

The notions of splitting field and generic splitting
field naturally extend to this setting.  That is, given
$\alpha \in \mathcal{F}(E)$, we will say that $L/E$ is
a {\em splitting field} for $\alpha$ if $\alpha_L = 1_L$. We will
call a splitting field $K/E$ for $\alpha$
{\em generic} (respectively, {\em very generic}) 
if every splitting field $L/E$ for $\alpha$ is 
a specialization (respectively, a rational specialization)
of $K/E$. Moreover, we can now define $\cd(\alpha)$ by
\[ \text{$\cd(\alpha) = \min \, \{ \trdeg_E(K) \; | \; K/E$
is a generic splitting field for $\alpha \, \}$}  \]
and $\cd(\mathcal{F})$ by
\[ \text{$\cd(\mathcal{F}) = \max \, \{ \cd(\alpha) \, | \,
E/k$ is a finitely generated extension,
$\alpha \in \mathcal{F}(E)$ $\}$.} \]
Proposition~\ref{prop.splitting4} says that if $G$ is a connected
linear
algebraic group and $\mathcal{F} = H^1(_-, G)$ then the above
definition of $\cd(\alpha)$ agrees with Definition~\ref{def1}.
Moreover, Definition~\ref{def.cd} tells us that for this
$\mathcal{F}$, $\cd(\mathcal{F}) = \cd(G)$.

Note that none of the above definitions require $k$ to
be algebraically closed. In particular, it now makes sense
to talk about the canonical dimension of an algebraic group defined
over a non-algebraically closed field. This opens up interesting
directions for future research, but we shall not pursue
them in this paper. Instead we will continue to assume
that $k$ is algebraically closed, and our main focus will
remain on the functors $H^1(_-, G)$.  However, even in this
(more limited but already very rich) context, we will take
advantage of the notion of canonical dimension for a functor
by considering certain subfunctors of $H^1(_-, G)$.

We also remark that it is a priori possible that for some
functors $\mathcal{F}$, some fields $E/k$ and
some $\alpha \in \mathcal{F}(E)$ there will not exist a generic
splitting field; if this happens, then, according to our definition,
$\cd(\alpha) = \cd(\mathcal{F}) = \infty$.  However,
Proposition~\ref{prop.splitting4} tells us that this
does not occur for any functor of the form $H^1(_-, G)$,
where $G$ is a linear algebraic group, and consequently,
for any of its subfunctors.

\begin{example} \label{ex.levi} Isomorphic functors
clearly have the same canonical dimension.
In particular, suppose $G$ is a linear algebraic group
and $U$ is a normal unipotent subgroup of $G$.
Then the natural map $H^1(_-, G) \lra H^1(_-, G/U)$ is 
an isomorphism (see, e.g.,~\cite[Lemma 1.13]{sansuc})
and hence, $\cd(G) = \cd(G/U)$.  Taking $U$ to be 
the unipotent radical of $G$, we see that
that $\cd(G) = \cd(G_{\operatorname{red}})$, where
$G_{\operatorname{red}}$ is the Levi subgroup of $G$.
\end{example}

The following simple lemma slightly extends the observation
that isomorphic functors have the same canonical dimension.
This lemma will turn out to be surprisingly useful in the sequel.

\begin{lem} \label{lem.iso1}
Suppose $\tau \colon \mathcal{F}_1 \lra \mathcal{F}_2$ is a morphism
of functors with trivial kernel. Then for every finitely generated
field extension $E/k$,

\smallskip
(a) $\cd(\alpha)=\cd(\tau(\alpha))$
for any $\alpha\in \mathcal{F}_1(E)$.

\smallskip
(b) $\cd(\mathcal{F}_1)\leq\cd(\mathcal{F}_2)$.

\smallskip
(c) Moreover, if $\tau$ is surjective then $\cd(\mathcal{F}_1)=
\cd(\mathcal{F}_2)$.
\end{lem}

\begin{proof}
Since $\tau$ has trivial kernel, $\alpha$ and $\tau(\alpha)$ have
the same splitting fields and hence, the same generic splitting
fields. This proves part (a). Parts (b) and (c) follow
from part (a) and the definition of $\cd(\mathcal{F})$.
\end{proof}

\begin{example} \label{ex.pso}
Recall that the cohomology set $H^1(_-, \PSO_{2n})$ classifies
pairs $(A, \sigma)$, where $A$ is a central simple
algebras of degree $2n$ with an orthogonal
involution $\sigma$ of determinant $1$; see~\cite[p. 405]{boi}.
(Note that~\cite{boi} uses the symbol $\PGO$ instead of $\PSO$.)
Consider the morphism of functors
$f \colon H^1(_-, \SO_{2n}) \lra H^1(_-, \PSO_{2n})$
sending a quadratic form
$q$ of dimension $2n$ to the pair $(\Mat_{2n}(K), \sigma_q)$,
where $\sigma_q$ is the involution of $\Mat_{2n}(K)$
associated to $q$. We claim that $f$ has trivial kernel.
Indeed, $q \in \Ker(K)$ $\Longleftrightarrow$ $q$ gives rise to the
standard (transposition) involution on $\Mat_{2n}(K)$ 
$\Longleftrightarrow$ $q$ is the $2n$-dimensional form $\langle a, a, \dots, a,
a \rangle$ for some $a \in K^*$; cf. e.g., \cite[p.~14]{boi}.
On the other hand, since we are assuming that $K$ contains 
an algebraically closed base field $k$ of characteristic zero 
(and in particular, $K$ contains a primitive $4$th root of unity),
the form $\langle a, a \rangle $ is hyperbolic
(cf., e.g., \cite[Theorem I.3.2]{lam}), and hence, so is $q$.
This shows that $f$ has trivial
kernel, as claimed. Lemma~\ref{lem.iso1}(b) now 
tells us that $\cd(\PSO_{2n}) \geq \cd(\SO_{2n})$.
\end{example}

\begin{example} \label{ex.spin_n}
The exact sequence
\[ 1 \lra \mu_2 \lra \Spin_n \stackrel{\pi}{\lra} \SO_n \lra 1 \]
of algebraic groups gives rise to the exact sequence
\[ \xymatrix{\SO_n(_-)
\ar[r]^{\delta} & H^1(_-,\mu_2)\ar[r]& H^1(_-,\Spin_n)\ar[r]^{\pi_*}&
H^1(_-,\SO_n)& } \]
of cohomology sets, where $\delta$ is the spinor norm; see,
e.g., \cite[p. 688]{garibaldi}. 
Since $-1$ is a square in $k$, the unit form is hyperbolic, 
hence $\delta$ is surjective and thus $\pi_*$ has trivial 
kernel.  On the other hand,
the image of $\pi_*$ consists of quadratic forms $q$ 
of discriminant $1$ such that
\[ q' = \begin{cases} \text{$q$, if $n$ is even,} \\
\text{$q \oplus \langle 1 \rangle$, if $n$ is odd} \end{cases} \]
has trivial Hasse-Witt invariant.
Thus $\cd(\Spin_n) = \cd(\HW_n)$, where
$\HW_n$ is the set of $n$-dimensional quadratic forms $q$
such that $q'$ has trivial discriminant and trivial Hasse-Witt
invariant.
\end{example}

\begin{example} \label{ex.pf} Define the functors $\Pf_r$ and $\SPf_r$
by $\Pf_r(E)$ = $r$-fold Pfister forms defined over $E$ and
$\SPf_r(E)$ = scaled $r$-fold Pfister forms defined over $E$.
In other words,
\[ \SPf_r(E) = \{ 
\langle c \rangle \otimes q \, | \, c \in E^*, \; q \in \Pf_r(E) \} \, . \]
Taking $c = 1$ above, we see that
$\Pf_r$ is a subfunctor of $\SPf_r$; hence,
$\cd(\Pf_r) \le \cd(\SPf_r)$. On the other hand, since
$q$ and $\langle c \rangle \otimes q$ have the same splitting fields
for every $q \in \Pf_r(E)$ and every $c \in E^*$, we actually have
equality $\cd(\Pf_r) = \cd(\SPf_r)$.

Now suppose $q \in \Pf_r(E)$. Let $q'$ be a subform of $q$
of dimension $2^{r-1} + 1$. 
The argument in \cite[p. 29]{ks} shows that $K = E(q')$
is a generic splitting field for $\alpha$. Recall that
$E(q')$ is defined as the function field of the
quadric hypersurface $q' = 0$ in $\bbP_E^{2^{r-1}}$;
in particular, $\trdeg_E \, E(q) = 2^{r-1} - 1$.
Proposition~\ref{prop.splitting4} now tells us that
$\cd(q) \le 2^{r-1} - 1$. On the other hand, if $q$ is
anisotropic, a theorem of Karpenko and
Merkurjev~\cite[Theorem 4.3]{km} tells us that, in
fact $\cd(q) = 2^{r-1} - 1$. We conclude that
\begin{equation} \label{e.pf}
\cd(\SPf_r) = \cd(\Pf_r) = 2^{r-1} - 1 \, .
\end{equation}

We remark that the setting considered by Karpenko 
and Merkurjev in~\cite{km} is a bit different 
from ours in that they call a field $K/E$ splitting
for a quadratic form $q/E$ if $q_K$ is isotropic, where as
we use this term to indicate that $q_K$ is hyperbolic.
However, if $q$ is a Pfister form then the two
definitions coincide; cf., e.g., \cite[Corollary 10.1.6]{lam}. 
\end{example}

\begin{example} \label{ex.g2}
Consider the exceptional group $G_2$.  The functors $H^1(_-, G_2)$ and
$\Pf_3$ are isomorphic; see, e.g., \cite[Corollary 33.20]{boi}.
Hence, $\cd(G_2) = \cd(\Pf_3) = 3$; see~\eqref{e.pf}. 
\end{example}

\section{Groups of type $A$}
\label{sect.a_{n-1}}

In this section we will study canonical dimensions of the groups
$\GL_n/\mu_d$ and $\SL_n/\mu_e$, where $e$ divides $n$.
We define the functor
\begin{equation} \label{e.Cne}
C_{n, e} \colon \Fields \lra \Sets^*
\end{equation}
by $C_{n, e}(E/k) = \{$isomorphism classes of central 
simple $E$-algebras of degree $n$ and exponent dividing $e \}$.
The marked element in $C_{n, e}(E/k)$ is the split
algebra $\Mat_n(E)$. Clearly, $C_{n, e}$ is a subfunctor
of $H^1(_-, \PGLn)$.

\begin{lem} \label{lem6.2}
Let $n$ and $d$ be positive integers and $e$ be their
greatest common divisor. Then
$\cd(\GL_n/\mu_d) = \cd(\GL_n/\mu_e) =
\cd(\SL_n/\mu_e) = \cd(C_{n, d}) = \cd(C_{n, e})$.
\end{lem}

\begin{proof}
By Lemma~\ref{lem.slnd}, there are surjective morphisms of functors
\[ \begin{array}{l}
H^1(_-, \GL_n/\mu_d) \lra C_{n, d} \, , \\
H^1(_-, \GL_n/\mu_e) \lra C_{n, e} \quad \text{and} \\
H^1(_-, \SL_n/\mu_e) \lra C_{n, e}
\end{array} \]
with trivial kernels. Basic properties of the index
and the exponent of a central simple algebra tell
us that $C_{n, d} = C_{n, e}$; the rest follows from
Lemma~\ref{lem.iso1}(c).
\end{proof}

\begin{lem} \label{lem6.2a}
Let $n$ and $e$ be positive integers such that $e$ divides $n$,

\smallskip
(a) If $e' \, | \, e$ and $n' \,  | \, n$ then
$\cd(\SL_n/\mu_e) \ge \cd(\SL_{n'}/\mu_{e'})$.

\smallskip
(b) Suppose $n = n_1 n_2$ and $e = e_1 e_2$, where
$e_i \, | \, n_i$ and $n_1, n_2$ are relatively prime.  Then
\[ \cd(\SL_{n}/\mu_{e}) = \cd(\SL_{n_1}/\mu_{e_1} \, \times \,
\SL_{n_2}/\mu_{e_2}) \le
\cd(\SL_{n_1}/\mu_{e_1}) + \cd(\SL_{n_2}/\mu_{e_2}) \, \]

\smallskip
(c) Let $n = \prod p^{a_p}$ be the prime factorization of $n$
(here the product is taken over all primes $p$ and $a_p = 0$ for all but
finitely many primes) and $m = \prod_{p | e} p^{a_p}$.  Then
$\cd(\SL_n/\mu_e) = \cd(\SL_m/\mu_e)$,
\end{lem}

\begin{proof} (a) The morphism of functors
$C_{n', e'} \lra C_{n, e}$ given by
$A \mapsto \M_{\frac{n}{n'}}(A)$ has trivial kernel.
By Lemma~\ref{lem.iso1}(b),
$\cd(C_{n, e}) \ge \cd(C_{n', e'})$. The desired
inequality now follows from Lemma~\ref{lem6.2}.

\smallskip
(b) First note that the functors
\[ \text{$H^1(_-, \SL_{n_1}/\mu_{e_1} \times \SL_{n_2}/\mu_{e_2})$ and
$H^1(_-, \SL_{n_1}/\mu_{e_1}) \times H^1(_-, \SL_{n_2}/\mu_{e_2})$} \]
are isomorphic.
Thus by Lemma~\ref{lem.slnd}, there is a surjective morphism of functors
\[ H^1(_-, \SL_{n_1}/\mu_{e_1} \times \SL_{n_2}/\mu_{e_2}) \lra
C_{n_1, e_1} \times C_{n_2, e_2} \]
with trivial kernel. By Lemma~\ref{lem.iso1}(c),
\[ \cd(\SL_{n_1}/\mu_{e_1} \times \SL_{n_2}/\mu_{e_2}) = \cd(C_{n_1, e_1}
\times C_{n_2, e_2}) \]
and by Lemma~\ref{lem6.2}, $\cd(\SL_{n}/\mu_{e}) = \cd(C_{n, e})$.
The equality in part (b) now follows from the fact that
for relatively prime $n_1$ and $n_2$ the functors
$C_{n_1, \, e_1} \times C_{n_2, \, e_2}$ and
$C_{n, e}$ are isomorphic via $(A_1, A_2) \mapsto A_1 \otimes A_2$.
The inequality in part (b) is a special case
of Lemma~\ref{lem6.1}(d). 

\smallskip
(c) By Lemma~\ref{lem6.2}, $\cd(\SL_n/\mu_e) = \cd(C_{n, e})$
and $\cd(\SL_m/\mu_e) = \cd(C_{m, e})$.
On the other hand, basic properties of the index and the exponent
of a central simple algebra tell us that
the functors $C_{m, e}$ and $C_{n, e}$ are isomorphic
via $A \mapsto \M_{n/m}(A)$.
\end{proof}

\begin{thm} \label{thm.sl}
Let $\alpha \in H^1(E, \PGLn)$ be the class of 
a division algebra $A$ of degree $n = p^i$, where $p$ is
a prime.  Then $\cd(\alpha) = n - 1$.
\end{thm}

Let $X$ be the Brauer-Severi variety of $A$. 
By a theorem of Karpenko~\cite[Theorem 2.1]{karpenko1}
every rational map $X \dasharrow X$ defined over $E$ is necessarily 
dominant. (For a related stronger result, 
see Merkurjev~\cite[Section 7.2]{merkurjev2}.)
Theorem~\ref{thm.sl} is an easy consequence 
of this fact; we outline the argument below.

\begin{proof} 
The function field $K = E(X)$ is a generic splitting field
for $A$; in particular, as we pointed out in Example~\ref{ex6.2a}, 
$\cd(\alpha) \le n-1$. To prove the opposite inequality, assume 
the contrary: $A$ has a generic splitting field
$L/E$ of transcendence degree $<  n-1$. Let $Y$ be 
a variety (defined over $E$) with function field $L$. 
Since $k(X)/E$ and $K/E$ are both generic splitting fields for $\alpha$, 
each is a specialization of the other.  Arguing as in 
the proof of Proposition~\ref{prop.splitting4}, 
we see that there exist rational
maps $f_1 \colon X \dasharrow Y$ and $f_2 \colon Y \dasharrow X$
such that the diagram
$$\xymatrix{X\ar@{-->}[r]^{f_1}\ar@{-->}[dr] & 
Y\ar@{-->}[r]^{f_2}\ar@{-->}[d]&X\ar@{-->}[dl]\\ & X/G & }$$
commutes. Moreover, after replacing $Y$ by the closure of 
the graph of $f_2$ in $Y \times X$, we may assume that $f_2$
is regular. Now $f_2 \circ f_1$ is a well defined rational map
$X \dasharrow X$, and
\[ \dim \, f_2 \circ f_1(X) \le \dim(Y) < n-1 = \dim(X) \, , \]
contradicting Karpenko's theorem. This concludes the proof 
of Theorem~\ref{thm.sl}.
\end{proof}

\begin{cor} \label{cor6.3} 
Suppose $n = p^i n_0$ and $e = p^j$, where $\gcd(p, n_0) = 1$ and
$i \ge j$. Then $\cd(\SL_n/\mu_e) = \begin{cases} \text{$0$, if $j = 0$,} \\
        \text{$p^i - 1$, if $j \ge 1$.} \end{cases}$  
\end{cor}

\begin{proof} If $j = 0$ then
$\SL_n/\mu_e = \SL_n$ is special and hence, has canonical dimension $0$.
Thus we only need to consider the case where $j \ge 1$.  
By Lemma~\ref{lem6.2a}(c), $\cd(\SL_n/\mu_e) = \cd(\SL_{p^i}/\mu_{p^j})$.
Thus we may also assume that $n_0 = 1$, i.e., $n = p^i$.  

By Lemma~\ref{lem6.2}, $\cd(\SL_n/\mu_e) = \cd(C_{n, e})$.
Since $C_{n, e}$ is, by definition, a subfunctor of $H^1(_-, \PGLn)$,
Example~\ref{ex6.2} tells us that $\cd(\SL_n/\mu_e) \le n-1 = p^i - 1$. 
To prove the opposite inequality, let $A$ be
a division algebra of degree $p^i$ and exponent $p^j$
(such algebras are known to exist; see, e.g., \cite[Appendix 7C]{rowen})
and let $\alpha$ be the class of $A$ in $H^1(E, \PGLn)$. Then
by Theorem~\ref{thm.sl}, 
$\cd(\SL_n / \mu_e) \ge \cd(\alpha) = n - 1$, as desired.
\end{proof}

\begin{cor} \label{cor.sympl} 
Let $n = 2^i n_0$, where $i \ge 1$ and $n_0$ is odd.
Then $\cd(\PSymp_n) = 2^i - 1$.
\end{cor}

\smallskip
Here $\PSymp$ stands for the projective symplectic group. Note 
that these groups are sometimes denoted by the symbol
$\operatorname{PGSp}$; see, e.g.,~\cite[p. 347]{boi}.

\smallskip
\begin{proof} Recall that every $\alpha \in H^1(_-,\PSymp_{n})$
is represented by a pair $(A,\sigma)$, where $A$ is a central
simple algebra of degree $2n$ and exponent $\le 2$,
and $\sigma$ is a symplectic involution on $A$.
A central simple algebra has a symplectic involution if
and only if its exponent is $1$ or $2$; moreover, a symplectic
involution of a split algebra is necessarily hyperbolic.
In other words, the morphism of functors
\[ H^1(_-, \PSymp_n) \;  \lra \; C_{n, 2} \, , \]
given by $\alpha \mapsto A$ is surjective and has trivial kernel.
Here $C_{n, 2}$ is the functor of central simple algebras of
degree $n$ and exponent dividing $2$, as in~\eqref{e.Cne}. Thus
\[ \cd(\PSymp_{n}) = \cd(C_{n, 2})
\stackrel{\text{\tiny by Lemma~\ref{lem6.2}}}{=} \cd(\SL_{n}/\mu_2) 
\stackrel{\text{\tiny by Corollary~\ref{cor6.3}}}{=} 2^i - 1 \, , \]
as claimed.
\end{proof}

\section{Orthogonal and Spin groups}
\label{sect.son}

\begin{lem} \label{lem.so0}
(a) $\cd(\SO_{n-1}) \le \cd(\SO_n)$ for every $n \ge 2$.
Moreover, equality holds if $n$ is even.

\smallskip
(b) $\cd(\Spin_{n-1}) \le \cd(\Spin_n)$ for every $n \ge 2$.
Moreover, equality holds if $n$ is even.

\smallskip
(c) $\cd(\SO_n) \ge \cd(\Spin_n)$ for every $n \ge 2$.
Moreover, if $n \ge 2^r$, where $r \ge 3$ is an integer, then
$\cd(\Spin_n) \ge 2^{r-1} - 1$.
\end{lem}

\begin{proof}
(a) The morphism $\tau \colon H^1(_-,\SO_{n-1})\to H^1(_-,\SO_n)$,
sending a quadratic form $q$ to $\langle 1 \rangle \oplus q$
has trivial kernel. Lemma~\ref{lem.iso1}(b) now tells us that
$\cd(\SO_{n-1}) \le \cd(\SO_n)$.

To prove the opposite inequality for
$n$ is even, let $q = \langle a_1, \dots, a_n \rangle \in
\SO_n$. Then $q = \langle a_1 \rangle \otimes \tilde{q}$, where
$\tilde{q} = \langle 1, a_1a_2,
\dots, a_1a_n \rangle$ lies in the image of $\tau$.
(Note that here we use the assumption that $n$ is even to conclude
that $\tilde{q}$ has discriminant $1$.)
Since $q$ and $\tilde{q}$ have the same splitting fields,
$\cd(q) = \cd(\tilde{q})$. On the other hand, since
$\tilde{q}$ lies in the image of $\tau$, Lemma~\ref{lem.iso1}(a)
tells us that $\cd(\tilde{q}) \le \cd(\SO_{n-1})$.
Thus $\cd(q) \le \cd(\SO_{n-1})$ and consequently,
$\cd(\SO_n) \le \cd(\SO_{n-1})$, as desired.

\smallskip
(b) is proved by the same argument as (a), using
the identity $\cd(\Spin_n) = \cd(\HW_n)$
of Example~\ref{ex.spin_n}.
The first assertion follows from the fact that
$\tau$ restricts to a morphism $\HW_{n-1} \lra \HW_n$.
In the proof of the second assertion, the key point is that
if a quadratic form $q = \langle a_1 \rangle \otimes \tilde{q}$,
of even dimension $n$, has trivial
discriminant and trivial Hasse-Witt invariant then so does
$\tilde{q}$; the rest of the argument goes through unchanged.

\smallskip
(c) The first inequality follows from the fact that
$\HW_n$ is a subfunctor of $H^1(_-, \SO_n)$. To prove
the second inequality, note that by part (b) we may assume
$n = 2^r$. Since the discriminant and the Hasse-Witt
invariant of an $r$-fold Pfister form are both trivial for
any $r \ge 3$, we see that $\Pf_r$ is a subfunctor of $\HW_n$
and thus
\[ \cd(\Spin_n) = \cd(\HW_n) \ge \cd(\Pf_r) = 2^{r-1} - 1 \, ; \]
see Example~\ref{ex.pf}.
\end{proof}

\begin{example} \label{ex.so1}
(a) $\cd(\SO_n) = \begin{cases} 
\text{0, if $n = 1$ or $2$,} \\
\text{1, if $n = 3$ or $4$,} \\
\text{3, if $n = 5$ or $6$.} \end{cases}$

\smallskip
(b) $\cd(\Spin_n) = 0$ for $n = 3$, $4$, $5$ or $6$.

\smallskip
(c) $\cd(\Spin_n) = 3$ for $n = 7$, $8$, $9$ or $10$.
\end{example}

\begin{proof} (a) Note that $\SO_1 = \{ 1 \}$, $\SO_3 \simeq
\PGL_2$, and $\SO_6 \simeq \SL_4/\mu_2$. Hence,  $\cd(\SO_1) = 0$, and  
(by Corollary~\ref{cor6.3}) $\cd(\SO_3) = 1$ and $\cd(\SO_6) = 3$.
The remaining cases follow from Lemma~\ref{lem.so0}(a).

\smallskip
(b) In view of Lemma~\ref{lem.so0}(b), 
it is enough to show that $\cd(\Spin_6) = 0$.
By the Arason-Pfister theorem~\cite[Theorem 10.3.1]{lam},
the only $6$-dimensional form
with trivial discriminant and trivial Hasse-Witt invariant is the split
form. In other words, $\HW_6$ is the trivial functor and thus
\[ \cd(\Spin_6) = \cd(\HW_6) = 0 \, . \]

Alternative proof of (b): Exceptional isomorphisms
of simply connected simple groups tell us that
$\Spin_3 \simeq \SL_2$, $\Spin_5 \simeq {\rm Sp}_4$ and
$\Spin_6 \simeq \SL_4$ are all special and hence, have
canonical dimension $0$; cf. Lemma~\ref{lem6.1}(c).

\smallskip
(c) Using the Arason-Pfister theorem once again, we see that
every $8$-dimensional quadratic form with trivial
discriminant and trivial Hasse-Witt invariant, is a
scaled Pfister form; see~\cite[Corollary 10.3.3]{lam}.
Thus
\[ \cd(\Spin_7) = \cd(\Spin_8) = \cd(\SPf_3) = 3 \, ; \]
see Example~\ref{ex.pf}. On the other hand, by a theorem
of Pfister every $q \in \HW_{10}$ is isotropic, i.e., has 
the form $\langle 1, - 1 \rangle \oplus q'$, where $q' \in \HW_8$;
see~\cite[Proof of Satz 14]{pfister} (cf. also~\cite[Theorem 4.4]{km}).

Applying Lemma~\ref{lem.iso1}(a) to the morphism $\HW_8 \lra \HW_{10}$
given by $q' \lra \langle 1, - 1 \rangle \oplus q'$, we see that
\[ \cd(q) = \cd(q') \le \cd(\HW_8) = \cd(\Spin_8) \, . \]
This shows that $\cd(\Spin_9) = \cd(\Spin_{10}) \le \cd(\Spin_8) = 3$.
The opposite inequality is given by Lemma~\ref{lem.so0}(b).
\end{proof}

\begin{prop} \label{prop.so}
$\cd(\SO_{2m}) \le \frac{m(m-1)}{2}$ for every $m \ge 1$.
\end{prop}

\begin{proof}
Write $\SO_{2m} =\SO(q)$, where $q$ is a non-degenerate
quadratic form on $k^{2m}$. Let $X = \operatorname{Gr}_{iso}(m, 2m)$
be the Grassmannian of maximal (i.e., $m$-dimensional) 
$q$-isotropic subspaces of $k^{2m}$, i.e., 
of $m$-dimensional subspaces contained in the quadric $Q \subset k^{2m}$
given by $q = 0$. It is well known that $X$ is a projective variety 
with two irreducible components $X_1$ and $X_2$, each 
of dimension $\frac{m(m-1)}{2}$; see e.g.,~\cite[Section 6.1]{gh}.
Using the Witt Extension Theorem (see, e.g.,~\cite[p. 26]{lam}),
it is easy to see that the full orthogonal group $\Orth(q)$ 
acts transitively on $X$ and $\SO(q)$ acts transitively on 
each component $X_i$ ($i = 1, 2$).
Fix an isotropic subspace $L \in X_1$ and let $P = \Stab_{\SO(q)}(L)$. 
By Lemma~\ref{lem6.1}(a), with $G = \SO(q)$ and $H = P$, we have
\begin{eqnarray*} 
\cd(\SO_{2m}) \le \cd(P) + \dim(\SO_{2m}) - \dim(P) = \\
\cd(P) + \dim(X_1) = \cd(P) + \frac{m(m-1)}{2} \, . \end{eqnarray*} 
It remains to show that $\cd(P) = 0$. We claim that
the Levi subgroup of $P$ is naturally isomorphic to $\GL(L) \simeq \GL_m$
via $f \colon P \lra \GL(L)$, where
$f(g) = g_{| \, L}$. Once this claim is 
established, Example~\ref{ex.levi} tells us that
$\cd(P) = \cd(\GL_m) = 0$. (The last equality follows from the fact that
$\GL_m$ is special.)

To prove the claim, note that by the Witt Extension Theorem,
$f$ is a surjective homomorphism.  It remains to show 
that $\Ker(f)$ is unipotent. 
Indeed, choose a basis $e_1, \dots, e_{2m}$ of $k^{2m}$ 
so that 
\[ q (x_1 e_1 + \dots + x_{2m} e_{2m}) = x_1 x_{m+1} + \dots + x_m x_{2m} \] 
and $L$ is the span of $e_1, \dots, e_m$. Then every 
$g \in \Ker(f)$ has the form 
\begin{equation} \label{e.Gram0}
g = \begin{pmatrix} I_m & A \\
                       O_m & B \end{pmatrix} \, , 
\end{equation}
for some $m \times m$-matrices $A$ and $B$. (Here
$O_m$ and $I_m$ are, respectively, the zero and 
the identity $m \times m$-matrices.) 
The condition that $g \in \Orth(q)$ translates into 
\begin{equation} \label{e.Gram}
g \begin{pmatrix} O_m & I_m \\
   I_m & O_m \end{pmatrix}  g^{\text{Transpose}} =
 \begin{pmatrix} O_m & I_m \\
   I_m & O_m \end{pmatrix} \, . 
\end{equation}
Substituting~\eqref{e.Gram0} into~\eqref{e.Gram}, 
we see that $B = I_m$. Formula~\eqref{e.Gram0} now shows that
$g$ is unipotent; consequently, $\Ker(f)$ is a unipotent group, 
as claimed.
\end{proof}

\begin{conjecture} \label{conj.so}
$\cd(\SO_{2m-1}) = \cd(\SO_{2m}) = \frac{m(m-1)}{2}$ for every $m \ge 1$.
\end{conjecture}

\begin{remark} \label{rem.karpenko} Conjecture~\ref{conj.so} was recently
proved by Karpenko~\cite{karpenko2}; for an alternative proof due 
to Vishik, see~\cite{vishik}.
\end{remark}

\section{Groups of low canonical dimension}
\label{sect.cd1}

\begin{thm}\label{cd1}
Assume that $G$ is simple. Then $\cd(G)=1$ if and only if $G \simeq
\SL_{2m}/\mu_2$ or $\PSymp_{2m}$, where $m$ is an odd integer.
\end{thm}

\begin{proof} First of all, observe that if 
$G \simeq \SL_{2m}/\mu_2$ or $\PSymp_{2m}$, with $m$ odd, 
then indeed, $\cd(G) = 1$; see Corollary~\ref{cor6.3} and
Corollary~\ref{cor.sympl}. Thus we only need to show that no
other simple group has this property.
Our proof relies on the classification of simple
algebraic groups; cf., e.g.,~\cite[\S 24 and 25]{boi}.
(Note that~\cite{boi} uses the symbols $\Orth^+$ and
$\operatorname{PGSp}$ instead of $\SO$ and $\PSymp$.
Recall also that we are working over an algebraically closed
base field $k$ of characteristic zero.)
We begin by observing
that $\cd(G) \ge 2$ for every simple group
of exceptional type; see Examples~\ref{ex.exceptional1},
\ref{ex.exceptional2} and~\ref{ex.g2}.

Now suppose $\cd(G) = 1$ and $G$ is of type $A$.
Then $G \simeq \SL_n/\mu_e$, where $e$ divides $n$.
Let $p$ be a prime dividing $e$.
Then we can write $e = p^j e_0$ and
$n = p^i n_0$, where $i \ge j \ge 1$, 
and $\gcd(p, e_0) = \gcd(p, n_0) = 1$. 
Then
\[ 1 = \cd(G) \stackrel{\text{\tiny by Lemma~\ref{lem6.2a}(a)}}{\ge}  
\SL_{p^i}/\mu_{p^j} 
\stackrel{\text{\tiny by Corollary~\ref{cor6.3}}}{\ge}  p^i - 1 \, , \]
which is only possible if $p = 2$ and $i = 1$. This 
implies that $e$ cannot be divisible by $4$ or by 
any prime $p \ge 3$; in other words, $e = 1$ or $2$. 
If $e = 1$ then $G = \SL_n$
is special and thus $\cd(G) = 0$. If $e = 2$ then 
$i = 1$ implies that $n \equiv 2 \, \pmod{4}$, as claimed.

Next suppose $G$ is of type $C$.
Then $G$ is isomorphic to $\Symp_{2m}$ or $\PSymp_{2m}$.
The groups $\Symp_{2m}$ are special and thus have 
canonical dimension $0$. By Corollary~\ref{cor.sympl}, 
$\cd(\PSymp_{2m}) = 1$
if and only if $m$ is odd. This completes the proof
of Theorem~\ref{cd1} for groups of type $A$ or $C$.

Now suppose $G$ is of type $B$ or $D$.
We have already considered some of these groups. 
In particular, 
\begin{itemize}

\smallskip
\item
$\cd(\SO_n) \ge \cd(\SO_5) = 3$ for any $n \ge 5$ (see Lemma~\ref{lem.so0}(a)
and Example~\ref{ex.so1}(a)),

\smallskip
\item
$\cd(\PSO_{2n}) \ge \cd(\SO_{2n}) \ge 3$
for any $n \ge 3$ (see Example~\ref{ex.pso}),

\smallskip
\item
$\cd(\Spin_n) = 0$ for $n = 3, 4, 5, 6$ (see Example~\ref{ex.so1}(b)),
and

\smallskip
\item
$\cd(\Spin_n) \ge \cd(\Spin_7) = 3$ for any $n \ge 7$
(see Lemma~\ref{lem.so0} and Example~\ref{ex.so1}(c)). 
\end{itemize}

We also remark that $\SO_2 \simeq \bbG_m$
and $\SO_4 \simeq (\SL_2 \times \SL_2)/\mu_2$ are not simple, and
$\SO_3 \simeq \PGL_2 = \SL_2/\mu_2$ was considered above.

This covers every simple group of type $B$; the only simple 
groups of type $D$ we have not yet considered
are $G = \Spin^{\pm}_{4n}$ ($n \ge 2$);
cf., e.g., \cite[Theorems 25.10 and 25.12]{boi}.
The natural projection $\pi \colon \Spin_{4n} \lra \SO_{4n}$ factors
through $\Spin_{4n}^{\pm}$:
\[ \pi \colon \Spin_{4n} \stackrel{f}{\lra}
\Spin_{4n}^{\pm} \lra \SO_{4n} \, . \]
Since $\pi_* \colon H^1(_-, \Spin_{4n}) \lra H^1(_-, \SO_{4n})$
has trivial kernel (see Example~\ref{ex.spin_n}),
so does $f_* \colon H^1(_-, \Spin_{4n}) \lra H^1(_-, \Spin_{4n}^{\pm})$.
Now for any $n \ge 2$,
\[ \cd(\Spin_{4n}^{\pm})
\stackrel{\text{\tiny by Lemma~\ref{lem.iso1}(b)}}{\ge}
\cd(\Spin_{4n})
\stackrel{\text{\tiny by Lemma~\ref{lem.so0}(b)}}{\ge}
\cd(\Spin_{8})
\stackrel{\text{\tiny by Example~\ref{ex.so1}(c)}}{=} 3 \, . \]
This completes the proof of Theorem~\ref{cd1}.
\end{proof}

\begin{remark} The above argument also shows that if a simple 
classical group $G$ has canonical dimension $2$ then either 
(i) $G \simeq \SL_{3m}/\mu_3$, where $m$ is prime to $3$
or possibly (ii) $G \simeq \SL_{6m}/\mu_6$, where $m$ is prime to $6$. 
In case (i), we know that $\cd(G) = \cd(\PGL_3) = 2$; 
see Corollary~\ref{cor6.3}.  In case (ii), $\cd(G) = \cd(\PGL_6)$ 
(see Corollary~\ref{lem6.2a}(c)); we do not know whether this number 
is $2$ or $3$.
\end{remark}

\section{The functor of orbits and homogeneous forms}
\label{sect.functor}

We now briefly recall the definition of essential dimension
of a functor, due to Merkurjev~\cite{merkurjev}.

Let $\mathcal{F}$ be a functor from the category
of all field extensions $K$ of $k$ to the category of sets.
(For our purposes, it is sufficient
to consider only finitely generated extensions $K/k$.)
Given $\alpha \in {\mathcal F}(K)$, we define $\ed(\alpha)$
as the minimal value of $\trdeg_k(K_0)$, where $k \subset K_0 \subset K$
and $\alpha$ lies in the image of the natural map ${\mathcal F}(K_0)
\lra {\mathcal F}(K)$. The essential dimension
$\ed(\mathcal{F})$ of the functor $\mathcal{F}$ is then defined
as the maximal value of $\ed(\alpha)$,
as $\alpha$ ranges over $\mathcal{F}(K)$ and $K$ ranges over
all field extensions of $k$.
In the special case, where $G$ is an algebraic group and
$\mathcal{F}= H^1(_-, G)$
we recover the numbers defined in Section~\ref{prel.ed}:
$\ed(\alpha) = \ed(X, G)$, where $\alpha \in H^1(K, G)$ and $X$ is
a generically free $G$-variety representing $\alpha$. Moreover,
$\ed(H^1(_-, G)) = \ed(G)$. For details, see~\cite{bf2}.

Now to each $G$-variety $X$ we will associate the functor
$\Orb_{X,G}$ given by $\Orb_{X,G}(L)=X(L)/ \sim$, where $a \sim b$
for $a, b \in X(L)$, if $a = g \cdot b$ for some $g \in G(L)$.
Given an $L$-point $a \in X(L)$, we shall denote $a
\pmod{\sim}$ by $[a] \in \Orb_{X, G}(L)$.
Using this terminology, Definition~\ref{def1} can be rewritten as follows.

\begin{prop}\label{cd_func} $\ed \, [\eta]=\cd(X,G)+\dim X/G$,
where $\eta \in X(k(X))$ is the generic point of $X$.
\end{prop}

\begin{proof} Let $Y$ be a variety with function field $k(Y) = L$.
Then $z \in X(L)$ may be viewed as a rational map $\phi_z \colon
 Y \dasharrow X$
and $g \in G(L)$ as a rational map $f_g \colon Y \dasharrow G$.
The point $g \cdot z$ of $X(L)$ corresponds to the map
$F_{z, g} \colon Y \dasharrow X$ given by $F_{z, g}(y) = f_g(y) \cdot 
\phi_z(y)$.
Consequently, the definition of $\ed \, [z]$ can be rewritten as
\begin{equation} \label{e.functorial}
\ed \, [z] = \min_{g \in G(L)} \, \{ \, \trdeg_k \, k(F_{z, g}(Y)) \, \} \, .
\end{equation}

Now set $z = \eta$, $L = k(X)$, $Y = X$, and $\phi = id_X$.
The element $g \in G(L)$ is then a rational map
$f = f_g \colon X \dasharrow G$,
$F = F_{z, g} \colon X \dasharrow X$ is, by definition,
a canonical form map (see Definition~\ref{def.fandF})
and the proposition follows from~\eqref{e.functorial} 
and Definition~\ref{def1}.
\end{proof}

\begin{remark} \label{rem.functor} 
By the definition of $\ed(\Orb_{X,G})$, we have
\[ \ed(\Orb_{X,G}) \ge \ed \, [\eta] \, . \]
We do not know whether or not equality holds in general. 
\end{remark}

For the rest of this paper we will focus on the following example.
Let $N = \begin{pmatrix} n + d -1 \\ d \end{pmatrix}$ and let
$X = \bbA^N$ be the space of degree $d$ forms
in $n$ variables $x = (x_1, \dots, x_n)$. That is, elements of
$\bbA^N$ are forms $p(x_1, \dots, x_n)$ of degree $d$ and
elements of $\bbP^{N-1}$ are hypersurfaces
$p(x_1, \dots, x_n) = 0$. The generic point of $\bbA^N$ is
the ``general" degree $d$ form in $n$ variables as
\[ \phi_{n, d}(x) =
\sum_{i_1 + \dots + i_n = d} a_{i_1, \dots, i_n} x_1^{i_1} \dots
x_n^{i_n} \in K[x_1, \dots, x_n] \, , \]
where $a_{i_1, \dots, i_n}$ are independent variables, $K = k(\bbA^N)$ 
is the field these variables generate over $k$.  Then
\begin{equation} \label{e.ed-phi}
\ed \, [\phi_{n, d}] = \min_{g \in \GL_n(K)}
\trdeg_k \, k(b_{i_1, \dots, i_d} \, | i_1 + \dots + i_n = d) \, ,
\end{equation}
where
\[ \phi_{n, d}(g \cdot x)
= \sum_{i_1 + \dots + i_n = d} b_{i_1, \dots, i_n} x_1^{i_1}
\dots x_n^{i_n} \, . \] 

The generic point of $\bbP^{N-1}$ is the ``general"
degree $d$ hypersurface $\phi_{n, d}(x) = 0$ in
$\bbP^{N-1}(K)$, which we denote by $H_{n, d}$. 

\begin{lem} \label{lem.hyper1a}
Let $N = \begin{pmatrix} n + d -1 \\ d \end{pmatrix}$ and
\[ D = \dim(\bbA^N/\GL_n) = \dim(\bbP^{N-1}/\PGLn) = 
\dim(\bbA^N/(\bbG_m \times \GL_n)) \, . \]
(Here $\bbG_m$ acts on $\bbA^N$ by scalar multiplication.) Then

\smallskip
(a) $\ed \, [\phi_{n, d}]  = D + \cd(\bbA^N, \GL_n)$.

\smallskip
(b) $\ed \, [H_{n, d}] =  D + \cd(\bbP^{N-1}, \GL_n) =   
 D + \cd(\bbP^{N-1}, \PGL_n) =$   

$\quad \quad \quad \quad \quad \;  
=   D + \cd(\bbA^{N}, \G_m \times \GL_n)$. 
\end{lem}

\begin{proof} Part (a) and the first equality in part (b)
are immediate consequences of Proposition~\ref{cd_func}.
To complete the proof of (b), note that
\[ \cd(\bbP^{N-1}, \GL_n) =  \cd(\bbP^{N-1}, \PGLn) \]
by Lemma~\ref{lem3}(b), and 
\[ \cd(\bbP^{N-1}, \GLn) = \cd(\bbA^{N}, \G_m \times \GL_n) \]
by Proposition~\ref{prop.product2} (with $H = \bbG_m \times \{ 1 \}$).
\end{proof}

\begin{cor} \label{cor.hyper1b}
Let $K = k(\bbA^N)$ (as above) and $K_0 = k(\bbP^N)$.
Given $g \in \GL(K)$,
let $F_g$ be the field extension of $k$ generated by elements
of the form 
\[ \text{$\frac{b_{i_1, \dots, i_n}}{b_{j_1, \dots, j_n}}$,  
where $\phi_{n, d}(g \cdot x)
= \sum_{i_1 + \dots + i_n = d} b_{i_1, \dots, i_n} x_1^{i_1}
\dots x_n^{i_n}$,} \]
and $i_1, \dots, i_n, j_1, \dots, j_n \ge 0$ satisfy
$i_1 + \dots + i_n = j_1 + \dots + j_n = d$ 
and $b_{j_1, \dots, j_n} \ne 0$.  Then
\[ \ed \, [H_{n, d}] = \min_{g \in \GL_n(K_0)} \, \trdeg_k (F_g) =
\min_{g \in \GL_n(K)} \trdeg_k (F_g) \, . \]
\end{cor}

\begin{proof} 
The first equality is an immediate consequence of
the definition of $\ed \, [H_{n, d}]$. 
To prove the second equality, we use the identity
$\ed \, [H_{n, d}] = D + \cd(\bbA^{N}, \G_m \times \GL_n)$ of
Lemma~\ref{lem.hyper1a}.  Proposition~\ref{cd_func} now tells us that
\[ \ed \, [H_{n, d}] = \min_{g \in \GL_n(K), \, c \in K^*}
\trdeg_k(c b_{i_1, \dots, i_n} \, | \, i_1 + \dots + i_n = d) \, . \]
The minimum is clearly attained if $c = b_{j_1, \dots, j_n}^{-1}$,
for some (and thus any) $j_1, \dots, j_n$ such that  
$b_{j_1, \dots, j_n} \ne 0$, and the corollary follows.
\end{proof}

In view of~\eqref{e.ed-phi}, 
it is natural to think of $\ed \, [\phi_{n, d}]$ 
as the minimal number 
of independent parameters required to define 
the general form of degree $d$ in $n$ variables. 
Corollary~\ref{cor.hyper1b} says that $\ed \, [H_{n, d}]$ can be similarly 
interpreted as the minimal number of independent parameters
required to define the general degree $d$ hypersurface
in $\bbP^{n-1}$. Comparing the expressions for these numbers 
given by~\eqref{e.ed-phi} and Corollary~\ref{cor.hyper1b},
we see that they are closely related.

\begin{cor} \label{cor.e.hyper1}
$\ed \, [H_{n, d}] \le \ed \, [\phi_{n, d}] \le \ed \, [H_{n, d}] + 1$. 
\qed
\end{cor}

\begin{remark} \label{rem.bf1} Lemma~\ref{lem.hyper1a}(b) shows
that the number $\ed \, [H_{n, d}]$ is the essential dimension 
of the generic form $\phi_{n, d}$ in the sense of~\cite{bf1},
i.e., the essential dimension of the $\bbG_m \times \GL_n$-orbit of
the generic point of $\bbA^N$. Note that the emphasis in~\cite{bf1}
is on the essential dimension of the functor 
${\bf Hypersurfaces}_{n, d} = \Orb_{\mathbb{A}^N, \bbG_m \times \GL_n}$
(which is denoted there by $\mathcal{F}_{d,n}$), 
and, more specifically, on the functor ${\bf Hypersurfaces}_{3, 3}$
(which is denoted there by ${\bf Cub}_3$). As we pointed out in
Remark~\ref{rem.functor},
$\ed({\bf Hypersurfaces}_{n, d}) \ge \ed \, [H_{n, d}]$ but we do not know 
whether or not equality holds.
\end{remark}

\section{Essential dimensions of homogeneous forms I}

\begin{thm} \label{thm.hyper}
Let $n$ and $d$ be positive integers such that
$d \ge 3$ and $(n, d) \ne (2, 3)$, $(2, 4)$ or $(3, 3)$.
Then
\[ \ed \, [H_{n, d}] =  N
- n^2 + \cd(\GL_n/\mu_d) =
N - n^2 + \cd(\SL_n/\mu_{\gcd(n, d)}) \, , \]
where $N = \begin{pmatrix} n + d -1 \\ d \end{pmatrix}$.
\end{thm}

\begin{proof} First observe that by Lemma~\ref{lem6.2},
$\cd(\GL_n/\mu_d) = \cd(\SL_n/\mu_{\gcd(n, d)})$, so
only the first equality needs to be proved.

Secondly, under our assumption on $n$ and $d$, the
$\PGL_n$-action on $\bbP^{N-1}$ is generically free.
For $n = 2$ this is classically known (cf., e.g.,~\cite[p. 231]{pv}),
for $n = 3$, this is proved in~\cite{baily} and for $n \ge 4$
in~\cite{mm}.  Substituting
\[ D = \dim(\bbP^{N-1}/\PGL_n) = \dim(\bbP^{N-1}) - \dim(\PGL_n)
= N - n^2 \]
into Lemma~\ref{lem.hyper1a}(b), we reduce the theorem
to the identity
\begin{equation} \label{e.thm.hyper}
\cd(\bbA^N, \bbG_m \times \GL_n) = \cd(\GL_n/\mu_d) \, .
\end{equation}
To prove~\eqref{e.thm.hyper}, observe that the normal subgroup 
\[ S = \{ (t^{-d}, t) \, | \, t \in \bbG_m \} 
\subset \bbG_m \times \bbG_m \subset \bbG_m \times \GL_n \]
acts trivially on $\bbA^N$. Since $S$ is special, we have
\begin{eqnarray} \label{e.hyper1c}
 \cd(\bbA^N, \bbG_m \times \GL_n)  
\stackrel{\text{\tiny by Lemma~\ref{lem3}(b)}}{=}  &  \\
\cd(\bbA^N, (\bbG_m \times \GL_n)/S) 
\stackrel{\text{\tiny by Definition~\ref{def.cd}}}{=} 
\cd((\bbG_m \times \GL_n)/S) & \nonumber \, ,  
\end{eqnarray}
where the last equality is a consequence of the fact that
the $\PGLn$-action on $\bbP^{N-1}$ (and hence, 
the $(\bbG_m \times \GL_n)/S$-action on $\bbA^N$) is 
generically free.

Finally, consider the homomorphism 
$\GL_n \lra (\bbG_m \times \GL_n)/S$
given by $g \mapsto (1, g)$, modulo $S$.
Since we are working over an algebraically closed 
field $k$ of characteristic zero, this homomorphism
is surjective, and its kernel is exactly $\mu_d$. 
Thus
$(\bbG_m \times \GL_n)/S \simeq \GL_n/\mu_d$.
Combining this with \eqref{e.hyper1c}, we obtain~\eqref{e.thm.hyper}.
\end{proof}

The results of Section~\ref{sect.a_{n-1}} can now be used to determine
$\ed \, [H_{n, d}]$ for many values of $n$ and $d$ (and produce estimates
for others). In particular, combining Theorem~\ref{thm.hyper} with
Corollary~\ref{cor6.3}, we deduce Theorem~\ref{thm.intro} stated 
in the Introduction. 

\smallskip
The number $\ed \, [\phi_{n, d}]$ appears to be harder to compute
than $\ed \, [H_{n, d}]$.  By Corollary~\ref{cor.e.hyper1},
$\ed \, [\phi_{n, d}] = \ed \, [H_{n, d}]$ or
$\ed \, [\phi_{n, d}] = \ed \, [H_{n, d}] + 1$, but
for general $n$ and $d$, we do not know which of these cases occurs.
One notable exception is the case where $n$ and $d$ are relatively 
prime.

\begin{cor} \label{cor.hyper1}
Suppose $d \ge 3$, $\gcd(n, d) = 1$
and $(n, d) \ne (2, 3)$. Then

\smallskip
(a) $\ed \, [H_{n, d}] =
\begin{pmatrix} n + d -1 \\ d \end{pmatrix} - n^2$ and

\smallskip
(b) $\ed \, [\phi_{n, d}] =
\begin{pmatrix} n + d -1 \\ d \end{pmatrix} - n^2 + 1$.
\end{cor}

\begin{proof} Part (a) is a special case of Theorem~\ref{thm.intro}
(with $j = 0$). We can also deduce it directly from
Theorem~\ref{thm.hyper} by noting that
$\SL_n$ is a special group and thus $\cd(\SL_n) = 0$.

\smallskip
(b) In view of Corollary~\ref{cor.e.hyper1}, we only need to prove that
$\ed \, [\phi_{n, d}] \ge \ed \, [H_{n, d}] + 1$ or equivalently,
$\cd(\bbA^N, \GL_n) \ge 1$; see Lemma~\ref{lem.hyper1a}(a).
Recall that the central subgroup $\mu_d$ of $\GL_n$ acts trivially
on $\bbA^N$, and (under our assumptions on $n$ and $d$)
the induced $\GL_n/\mu_d$-action is generically free.
Thus the stabilizer in general position for the $\GL_n$-action
on $\bbA^N$ is $\mu_d$, and 
\[ \cd(\bbA^N, \GL_n) 
\stackrel{\text{\tiny by Proposition~\ref{prop.lower1}(c)}}{\ge} 
\ed(\mu_d) 
\stackrel{\text{\tiny by \cite[Theorem 6.2]{br1}}}{=} 1 \, , \]
as claimed.
\end{proof}

\begin{remark} \label{rem.product2} 
We have proved that if $n$ and $d$ are relatively 
prime then $\cd(\bbP^{N-1}, \GL_n) = 0$ but
$\cd(\bbA^N, \GL_n) = 1$. In particular, this shows that
the equality $\cd(X, G) = \cd(X/H, G)$ of Proposition~\ref{prop.product2}
fails for $X = \bbA^N$, $G = \GL_n$ and $H = \bbG_m$. Note that
Proposition~\ref{prop.product2} does not apply in this situation
because the $H$-action on $X$ is not generically free
(the subgroup $\mu_d$ acts trivially).
\end{remark}

\section{Essential dimensions of homogeneous forms II}
\label{sect.hom2}

In this section we will study $\ed \, [\phi_{n, d}]$ and $\ed \, [H_{n, d}]$
for the pairs $(n, d)$ not covered by Theorem~\ref{thm.hyper}.
We begin with a simple lemma.

\begin{lem} \label{lem.binary}
$\ed \, [H_{2, d}] \le d - 2$ for any $d \ge 3$.
\end{lem}

In the sequel we will only need this lemma for $d = 3$ and $4$.
For $n \ge 5$, Theorem~\ref{thm.intro} (with $n = 2$) gives a
stronger result, namely,
\[ \ed \, [H_{2, d}] = \begin{cases} \text{$d - 2$ if $d$ is even}, \\ 
\text{$d-3$ if $d$ is odd.} \end{cases} \]
However, the proof of the lemma below is valid for any $d \ge 3$.

\begin{proof}
The linear transformation $g \in \GL_2(K)$ given by
\[ x_1 \mapsto x_1 - \frac{a_{d-1, 1}}{n \, a_{d, 0}} \, , \quad
x_2 \mapsto x_2 \]
reduces the generic binary form
\[ \phi_{2, d}(x_1, x_2) = a_{d, 0} x_1^d + a_{d-1, 1} x_1^{d-1} x_2 +
\dots + a_{0, d} x_2^d \]
to 
\[ \phi_{2, d}(g \cdot (x_1, x_2)) = 
b_{d, 0} x_1^d + b_{d - 2, 2} x_1^{d-2} x_2^2 + \dots + 
b_{1, d-1} x_1 x_2^{d-1} + b_{0, d} x_2^d  \]
for some $b_{i, d-i} \in K = k(a_{0, d}, \dots, a_{d, 0})$.
Composing this linear transformation with
\[ x_1 \mapsto \frac{b_{0, d}}{b_{1, d-1}} x_1 \, , \quad
x_2 \mapsto x_2 \]
(and, by abuse of notation, denoting the composition by $g$ once again),
we may further assume $b_{1, d-1} = b_{0, d}$. 
The field $F_g = k(b_{i, d-i}/b_{0, d} \, | \, i = 1, \dots, d)$,
defined in the statement of Corollary~\ref{cor.hyper1b}, 
now has transcendence degree $\le d-2$. By
Corollary~\ref{cor.hyper1b} we conclude 
that $\ed \, [H_{2, d}] \le n - 2$.
\end{proof}

We are now ready to proceed with the main result of this section.

\begin{prop} \label{prop.hyper}
(a) $\ed \, [\phi_{n, 1}] = \ed \, [H_{n, 1}] = 0$.

\smallskip
(b) $\ed \, [\phi_{n, 2}] = n$ and $\ed \, [H_{n, 2}] = n-1$.

\smallskip
(c) $\ed \, [\phi_{2, 3}] = 2$ and $\ed \, [H_{2, 3}] = 1$.

\smallskip
(d) $\ed \, [\phi_{2, 4}] = 3$ and $\ed \, [H_{2, 4}] = 2$.

\smallskip
(e) $\ed \, [H_{3, 3}] = 3$.
\end{prop}

\smallskip
\noindent
We do not know whether $\ed \, [\phi_{3, 3}]$ is $3$ or $4$.

\begin{proof}
(a) A linear form $l(x_1, \dots, x_n)$ over $K$ can
be reduced to just $x_1$ by applying a linear transformation
$g \in \GL_n(K)$. Thus $\ed \, [\phi_{n, 1}] = \ed \, [H_{n, 1}] = 0$.

\smallskip
(b) Here $d = 2$, $N = n(n+1)/2$, and elements of $\bbA^N$ are
quadratic forms in $n$ variables. Diagonalizing the generic
quadratic form $\phi_{n, 2}$ over $K$, we see that
$\ed \, [\phi_{n, 2}] \le n$ and $\ed \, [H_{n, 2}] \le n - 1$.
In view of Corollary~\eqref{cor.e.hyper1} it suffices to show
that $\ed \, [\phi_{n, 2}] = n$.

The $\GL_n$-action on $\bbA^N$ has a dense orbit,
consisting of non-singular forms. In particular,
$D = \dim(\bbA^N/\GL_n) = 0$, so that by Lemma~\ref{lem.hyper1a}(a)
\[ \ed \, [\phi_{n, 2}] = \cd(\bbA^N, \GL_n) \, . \]
Since the stabilizer of a non-singular form is
the orthogonal group $\Orth_n$, $\bbA^N$ is birationally
$G$-equivariantly isomorphic to $\GL_n/\Orth_n$.
Thus
\begin{eqnarray*}
\ed \, [\phi_{n, 2}] = \cd(\bbA^N, \GL_n) = \cd(\GL_n / \Orth_n, \GL_n)
\stackrel{\text{\tiny by Corollary~\ref{cor4.3}(b)}}{=} \\
\ed(\Orth_n)
\stackrel{\text{\tiny by~\cite[Theorem 10.3]{reichstein}}}{=} n \, .
\end{eqnarray*}
This completes the proof of part (b).

\smallskip
(c) By Lemma~\ref{lem.binary}, $\ed \, [H_{2, 3}] \le 1$. Thus in view
of Corollary~\ref{cor.e.hyper1}, we only need to show that 
$\ed \, [\phi_{2, 3}] \ge 2$.

Here $N = 4$, and the $\GL_2$-action on $\bbA^4$
has a dense orbit consisting of
binary cubic forms with three distinct roots.
Applying Lemma~\ref{lem.hyper1a}(a), with $D = \dim (\bbA^4/\GL_2) = 0$,
as in part (b), we obtain
\[ \ed \, [\phi_{2, 3}] = \cd(\bbA^4, \GL_2) = \cd(\GL_2/S, \GL_2)
\stackrel{\text{\tiny by Corollary~\ref{cor4.3}(b)}}{=} 
\ed(S) \, , \]
where $S \subset \GL_2$ is the stabilizer of a binary cubic form
with three roots, say of $x^3 + y^3$.
Note that $S$ is a finite group and that matrices that
multiply $x$ and $y$ by third roots of unity form a subgroup of
$S$ isomorphic to $(\bbZ/ 3\bbZ)^2$.
Thus $\ed(S) \ge \ed \, (\bbZ/ 3\bbZ)^2 = 2$ (cf.
\cite[Lemma 4.1(a) and Theorem 6.1]{br1}), as desired.

\smallskip
(d) By Lemma~\ref{lem.binary}, $\ed \, [H_{2, 4}] \le 2$.
In view of Corollary~\ref{cor.e.hyper1}, it remains to prove the inequality
$\ed \, [\phi_{2, 4}] \ge 3$.
Note that since the invariant field $k(\bbA^5)^{\GL_2}$ is
generated by one element (namely, the cross-ratio
of the four roots of the quartic binary form), we
have $D = \dim(\bbA^5/\GL_2) = 1$. Thus we only need to show that
\begin{equation} \label{e.binary4}
\cd(\bbA^5, \GL_2) \ge 2 \, .
\end{equation}
Let $S$ be the stabilizer of
$f \in \bbA^5$ (i.e., of a degree 4 binary form) in general
position. By Proposition~\ref{prop.lower1}(c),
\[ \cd(\bbA^5, \GL_2) \ge \ed(S) \, . \]
To compute $\ed(S)$, recall that the stabilizer of $f$ in $\PGL_2$
is isomorphic to $(\bbZ/2 \bbZ)^2$; cf.~e.g.,~\cite[p. 231]{pv}.
It is now easy to see that $S$ fits into the sequence
\[ \{ 1 \} \lra \bbZ/ 4 \bbZ \lra S \lra (\bbZ/ 2 \bbZ)^2
\lra \{ 1 \} \, . \]
In particular, $S$ is a finite group which admits a surjective homomorphism
onto $(\bbZ/2 \bbZ)^2$. Thus $S$ is neither cyclic nor odd dihedral and
consequently, $\ed(S) \ge 2$; see~\cite[Theorem 6.2(a)]{br1}.
This concludes the proof of~\eqref{e.binary4} and thus of part (d).
For the sake of completeness, we remark that since $S$ is a finite
subgroup of $\GL_2$, we also have $\ed(S) \le 2$ and thus $\ed(S) = 2$.

\smallskip
(e) Here $N = 10$, and the rational quotient $\bbP^9/ \GL_3$
is the $j$-line, so that $D = \dim \, \bbP^9/\GL_3 = 1$.
Thus we only need to show
\begin{equation} \label{e.binary5}
\cd(\bbP^9, \GL_3) = 2 \, .
\end{equation}
An element of $\bbP^9$ (i.e., a plane cubic curve) in general position
can be written as $F_{\lambda} = x^3 + y^3 + z^3 + 3 \lambda xyz$.
Denote the stabilizer of $F_{\lambda}$ by $S \subset \GL_3$.
We will deduce~\eqref{e.binary5} from Corollary~\ref{cor5.3}(b).
Indeed, let $N$ be the normalizer of $S$ in $\GL_3$.
Since $\GL_3$ is a special group, $e(\GL_3, S) = \ed(S)$,
$e(\GL_3, N) = \ed(N)$ (see Lemma~\ref{lem.e1}(c)), and
Corollary~\ref{cor5.3}(b) assumes the following form:
\[ \ed(S) \le \cd(\bbP^9, \GL_3) \le \ed(N) - \dim(S) + \dim(N) \, . \]
Let $\overline{S}$ and $\overline{N}$ be the images of $S$ and $N$
in $\PGL_3$, under the natural projection $\GL_3 \lra \PGL_3$.
Note that $\overline{S}$ is a finite group (this follows from the fact
that $D = \dim \, \bbP^9/\PGL_3 = 1$).  In particular, $\dim(S) = 1$.
It thus suffices to show:

\smallskip
(e$_1$) $\ed(S) \ge 2$,

\smallskip
(e$_2$) $\overline{N}$ is a finite subgroup of $\PGL_3$
(and consequently, $\dim (N) =1$). 

\smallskip
(e$_3$) $\ed(N) \le 2$.

\smallskip
The inequality
(e$_1$) is a consequence of~\cite[Corollary 7.3]{ry1}, with $G = S$ and
\begin{equation} \label{e.grp.H}
H = \mathopen< \diag ( 1, \zeta, \zeta^2), \sigma \mathclose>
\simeq (\bbZ/3 \bbZ)^2 \, , \end{equation}
where $\sigma$ is a cyclic permutation of the variables $x, y, z$, and
$\zeta$ is a primitive third root of unity.
(Note that \cite[Corollary 7.3]{ry1}
applies because $S$ has no non-trivial unipotent elements, and the
centralizer of $H$ in $S$ is finite.)

To prove (e$_2$), note that $\overline{N}$ is the normalizer
of $\overline{S}$ in $\PGL_n$. The natural 3-dimensional
representation of $S \subset \GL_3$ is irreducible (to see
this, restrict to the subgroup $H$ of $S$ defined 
in~\eqref{e.grp.H}).  Hence, by Schur's lemma, 
the centralizer $C_{\PGL_n}(\overline{S}) = \{ 1 \}$, so that
$\overline{N} = N_{\PGL_n}(\overline{S})/ C_{\PGL_n}(\overline{S})$.
The last group is naturally isomorphic to a subgroup of $\Aut(\overline{S})$,
which is a finite group. This proves (e$_2$).

To prove (e$_3$), consider the natural representation of
$N \subset \GL_3$ on $\bbA^3$. ($e_2$) implies that
this representation is generically free;
cf.~\cite[Section 1]{bf1}.  Consequently,
$\ed(N) \le 3 - \dim(N) = 2$.

This completes the proof of part (e).
\end{proof}

\end{document}